\documentclass[12pt]{amsart}
\usepackage{amsmath}
\usepackage{amsthm}
\usepackage{amsfonts}
\usepackage{hyperref}
\usepackage{epsfig}

\newtheorem{theorem}{Theorem}[section]
\newtheorem{proposition}{Proposition}[section]

\newtheorem{remark}{Remark}[section]
\newtheorem{definition}{Definition}[section]
\newtheorem{question}{Question}[section]
\newtheorem{conjecture}{Conjecture}[section]

\renewcommand{\epsilon}{\varepsilon}
\renewcommand{\thesection}{\arabic{section}}

\newcommand{\eqnsection}{
\renewcommand{\theequation}{\thesection.\arabic{equation}}
   \makeatletter
   \csname  @addtoreset\endcsname{equation}{section}
   \makeatother}
\eqnsection

\newcommand{\abs}[1]{\left\vert #1\right\vert}
\newcommand{\R}{\mathbb{R}}
\newcommand{\N}{\mathbb{N}}
\newcommand{\Z}{\mathbb{Z}} 
\newcommand{\PR}{\mathbb{P}}
\newcommand{\ES}{\mathbb{E}}
\newcommand{\LR}{\mathcal{L}}
\newcommand{\IG}{\mathfrak{I}}

\newcommand{\1}[1]{{\mathbf 1}{\{#1\}}}
\newcommand{\hatxi}{\zeta}
\newcommand{\backtrack}{\mathcal{B}\mathcal{K}}

\author[G. Ben Arous]{G\'erard BEN AROUS}
\address{G\'erard Ben Arous, Courant Institute of Mathematical Sciences, New York University, New York, New York 10012
USA} \email{benarous@cims.nyu.edu}

\author[A. Fribergh]{Alexander FRIBERGH}
\address{Alexander Fribergh, CNRS and Universit\'e de Toulouse, Institut de Math\'ematiques (CNRS UMR 5219),
31062, Toulouse, France} \email{alexander.fribergh@math.univ-toulouse.fr}


\begin{document}

\title[]{Biased random walks on random graphs}

\maketitle

\tableofcontents

\section{Introduction}

These notes cover one of the topics programmed for the St Petersburg School in Probability and Statistical Physics of June 2012. The aim is to review recent mathematical developments in the field of random walks in random environment (RWRE). For a detailed background on RWREs we refer the reader to~\cite{Zeitouni} or~\cite{SZ2}.

 Our main focus will be on directionally transient and reversible random walks on different types of underlying graph structures, such as $\Z$, trees and $\Z^d$ for $d\geq2$. 

Rather than speaking abstractly about the current state of the field, we decide to first dive into the heart of the subject by presenting rapidly a simple model which encapsulates most of the key questions we want to address in those notes. We feel this will give the reader a clear framework and provide an early motivation to understand the issues at hand.

\subsection{The key issues in one simulation}\label{sect_intro_biased}

We introduce informally the model of the biased random walk on a Galton-Watson tree with leaves (a more detailed description of this model can be found in Section~\ref{sect_arbre_model}).

Given a tree $T$ and $\beta>0$, we can define the $\beta$-biased random walk on the tree $T$. It is a Markov chain $(X_n)_{n\in \N}$ on the vertices of $T$, such that if $u\neq \text{root}$ and $u$ has $k$ children $v_1,\ldots, v_k$ and parent $\overleftarrow{u}$, then
\begin{enumerate}
\item $P^{T}[X_{n+1}=\overleftarrow{u}|X_n=u]=\frac 1 {1+ \beta k}$, 
\item $P^{T}[X_{n+1}= v_i |X_n=u] = \frac {\beta}{1+\beta k}$, for $1\leq i\leq k$,
\end{enumerate}
and from the root all transitions to its children are equally likely. The walk is typically started at the root.

\begin{figure}
\centering 
\epsfig{file=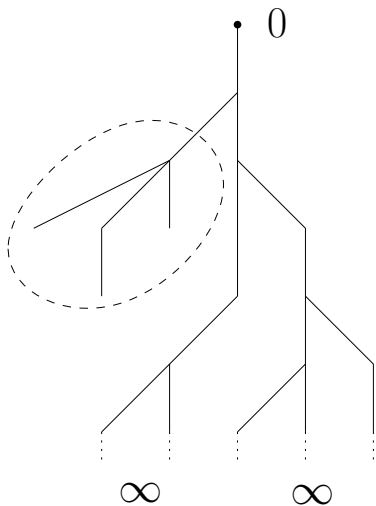, scale=0.8}
\caption{A Galton-Watson tree with the leaves.}
\end{figure}

The biased random walk on a Galton-Watson tree with leaves is obtained by choosing $T=T(\omega)$ to be a supercritical Galton-Watson tree with leaves (i.e.~positive probability of having no offspring) conditioned to be infinite. We use $p_k$ the probability to have $k$ offspring for the unconditioned Galton-Watson tree (in particular, our assumptions is that $p_0>0$).

We recall some of the most basic facts known about this walk (we will give details on this in Section~\ref{sect_arbre_model}). Firstly, there exists a constant $\beta_0>0$ such that
\[
\text{if $\beta> \beta_0$,} \qquad \lim_{n \to \infty}  d(\text{root},X_n) =\infty,
\]
by results in~\cite{lycap} and where $d(\cdot,\cdot)$ denotes the usual distance on trees. Secondly, for $\beta>\beta_0$, we have that (see~\cite{LPP})
\[
\lim_{n\to \infty}  \frac{d(\text{root},X_n)}n=v(\beta),
\]
where the limiting speed depends not only on $\beta$ but also on the offspring distribution. We choose to omit this dependence during the introduction.

The first result means that the walk is transient and the second means it has an asymptotic velocity $v(\beta)$ depending on the law of the tree and the bias.  Both results hold almost surely in the environment and the walk.

An interesting picture arises when one looks at simulations of the velocity (also called speed) $v(\beta)$ as a function of the bias for three different choices of trees, see Figure 2, 3, 4. In these simulations, the speed is estimated after 10 million steps.

\begin{remark}
It is important to emphasize that the speed is very hard to estimate through simulations for \lq\lq large\rq\rq\ values of the bias. Indeed, because of the trapping naturally occurring in this model, the convergence of $d(\text{root},X_n)/n$ towards the limiting velocity is very slow. The speed obtained through simulations is generally noticeably superior to the true speed for large values of the bias.

Nevertheless the simulations highlight several interesting features of this model and provide a very visual motivation for studying this model so we decided to include them. A qualitative representation of what the asymptotic speed should look like will be given later in these notes (see Figure~13).
\end{remark}

\begin{figure}
\centering 
\epsfig{file=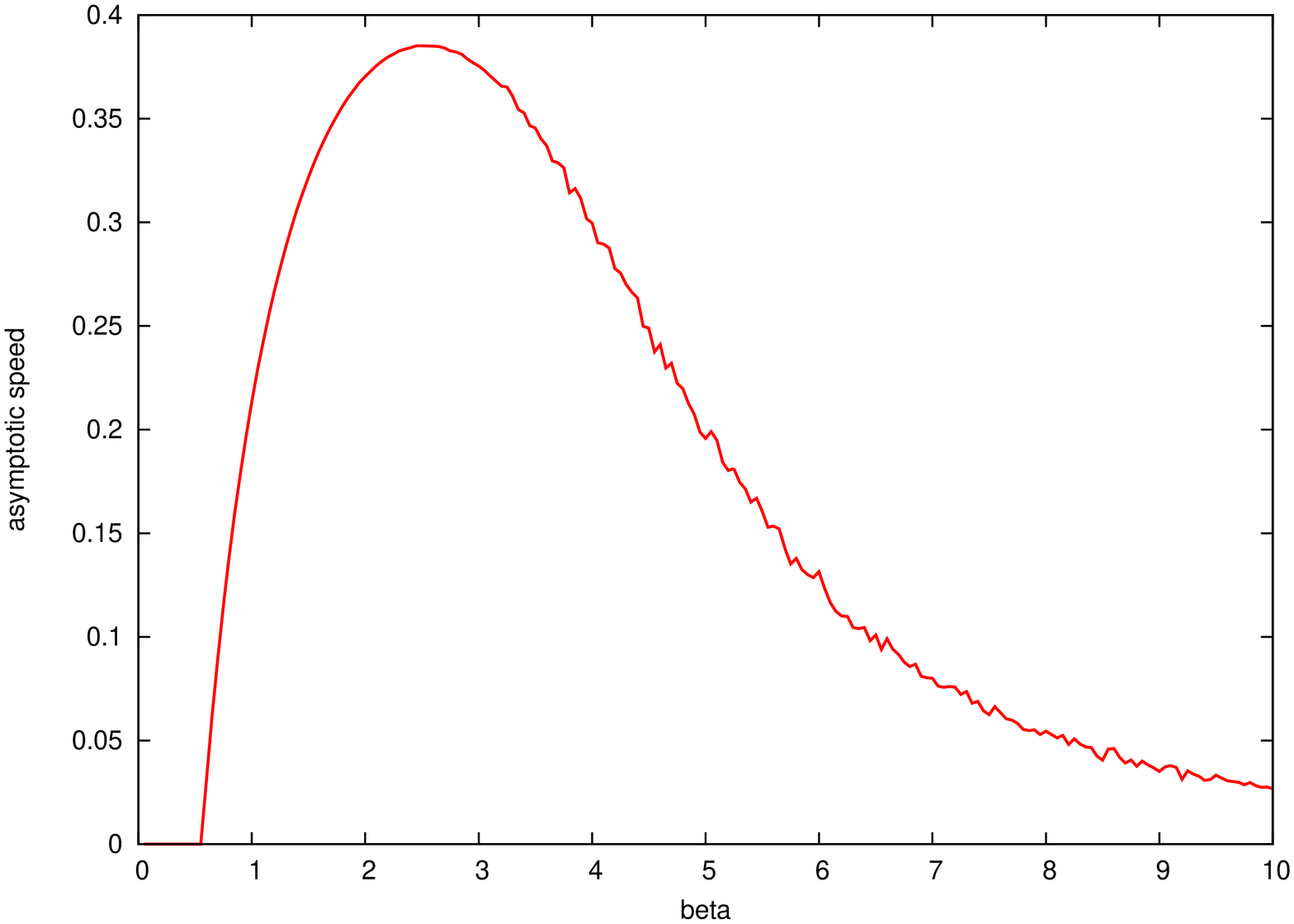, scale=0.5, angle=270}
\caption{Simulation of the speed of a biased random walk on a Galton-Watson tree with  leaves (by B. Rehle). Here $p_0=1/10$ and $p_2=9/10$.}

\epsfig{file=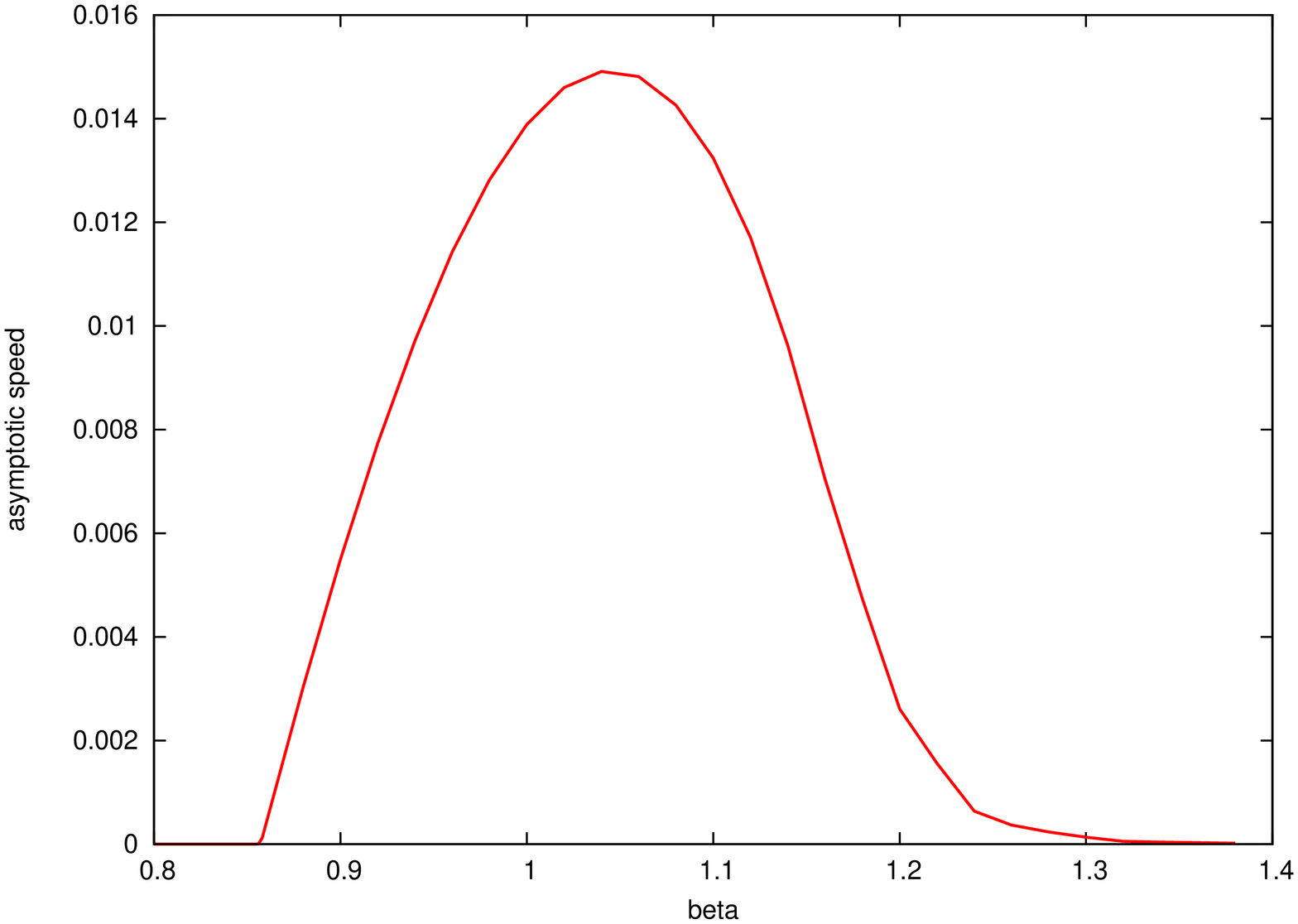, scale=0.5, angle=270}
\caption{Simulation of the speed of a biased random walk on a Galton-Watson tree with leaves (by B. Rehle). Here $p_0=1/4$, $p_1=1/3$ and $p_2=5/12$.}
\end{figure}

\begin{figure}
\centering 
\epsfig{file=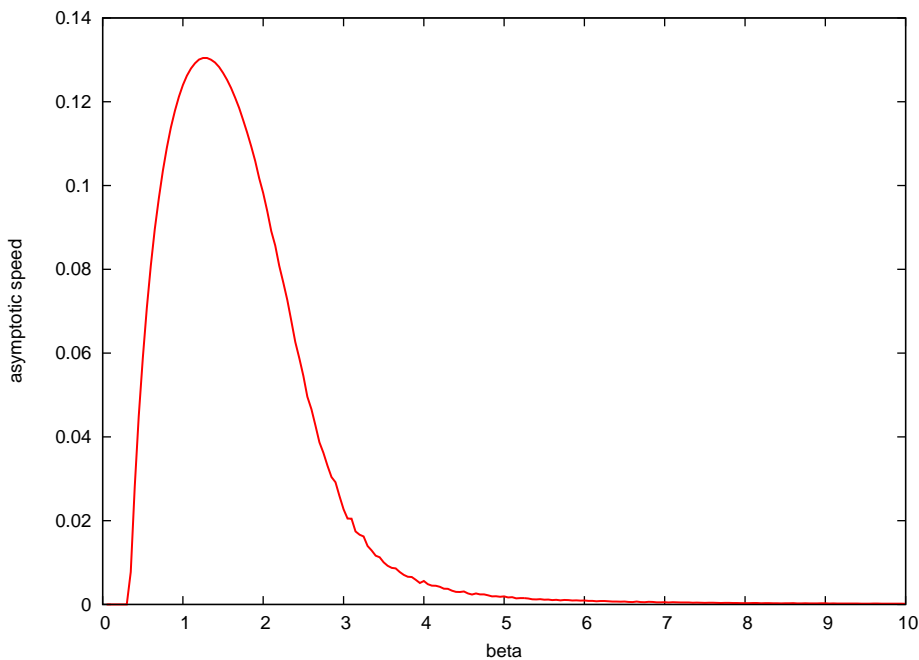, scale=0.5, angle=270}
\caption{Simulation of the speed of a biased random walk on a Galton-Watson tree with  leaves (by B. Rehle). Here $p_0=1/3$, $p_1=1/3$ and $p_8=1/3$}
\end{figure}

From these simulations, we seem to detect some similar patterns
\begin{enumerate} 
\item the speed first increases and then decreases (this is called unimodality),
\item the speed is eventually zero.
\end{enumerate}

At any point of the tree an increase in the bias will increase the local drift away from the root  (only at leaves does the drift remain unchanged). Yet, the walk is not necessarily sped up, actually, for stronger biases the motion eventually has zero speed.

The reason for the sub-ballistic regime is that the walk loses time in `traps' within the tree, from where it cannot go to infinity without having to battle for a long time against the drift which pushes it into the trap. The hypothesis $p_0>0$ is crucial for this to happen, since it is the leaves that create those dead-ends (see Figure 1). The case $\beta=\infty$ illustrates this phenomenon to the extreme since it turns the walk recurrent.

From these observations a set of very natural question appears,
\begin{enumerate}
\item Is it true that the speed is unimodal?
\item  Is it true that the speed is eventually zero?
\item How fast does the walker move when its speed is zero? After proper rescaling, what are the scaling limits of distance of the random walk from the origin?
\end{enumerate}

\subsection{Goal of the notes}\label{sect_plan}

The example we chose is convenient because it is simple to define mathematically and yet raises interesting questions. Nevertheless, the reader might question its relevance. As it turns out this model was considered as an approach to the more difficult problem of biased random walks on percolation clusters (a model presented in Section~\ref{perco}). This model already attracted the attention of physicists in the 80s, see~\cite{BD},~\cite{Dhar} and~\cite{DS}, not only for its interesting phenomenology but also because it presents similarities to concrete physical systems such as diffusion of particles in gels under gravity, or in centrifugal forces as in chromatographic columns, and hopping electron conduction in doped semiconductors in the presence of strong electric field. For a physical perspective on biased diffusion in disordered media, we invite the reader to take a look at Chapter 6 of the overview~\cite{Havlin}.

From  a more mathematical perspective, biased random walks on percolation clusters have also received a lot of attention.  This has been part of the larger project of trying to understand RWREs on $\Z^d$ (for $d\geq 2$), subject that has been very active in the last decade, following namely an important work of Sznitman and Zerner, see~\cite{SZ}. Although the field of RWREs is very wide, see~\cite{Zeitouni} or~\cite{SZ2} for overviews, we will chose to ignore a large part of it and focus on the case of reversible and directionally transient RWREs.

 These types of models have long been known to be related to trapping. But it is only in the last few years that concrete links have been established  between RWREs (that are directionally transient and reversible) and dynamics of traps models, more specifically a directed version of the  so-called Bouchaud Trap Model (BTM) (introduced in~\cite{Bouchaud}). This model, originally introduced to understand dynamics of spin glasses, turned out to be a key tool in understanding RWREs.

Our main goal is not to give a thorough overview of trap models, a topic covered in details in~\cite{BenarousCerny}. Rather, one of our goals will be to explain how a variety of models: biased random walks on Galton-Watson trees, on percolation clusters and many RWREs on $\Z$ can all be linked to one single trap model. This link is the crucial element to prove of a myriad of interesting properties such as anomalous scaling limits and aging for a wide class of RWREs. This will provide a very complete answer to the last two questions asked at the end of the previous section.

Another goal will be to study trapping under a different light for which RWREs offer a natural setting: understanding the influence of traps on the limiting velocity of random dynamics. In our concrete example of biased random walks on Galton-Watson tree with leaves, we can effectively tune the strength of the trapping by changing the value of $\beta$, and, as we saw in Figure 2-3-4, this can be felt in the behavior of the limiting velocity $v$ in function of $\beta$. The exact nature of this link is an interesting subject where deceivingly simple questions remain open.

\subsection{Plan of the notes}\label{sect_plan1}

We will take a progressive approach where we start from easier models and build up to more sophisticated examples. We start in Section~\ref{sect_1d} by presenting models where the underlying graph is one-dimensional, this will include a key toy model (called the directed Bouchaud Trap Model) and one-dimensional random walks in random environments. This first part is somewhat at odds with the rest of these notes that will focus on RWREs on trees and $\Z^d$. Nevertheless, it allows for a simple way to showcase the parallel between trap models and RWREs.

The second main section (Section~\ref{sect_arbre_model}) of this paper will be devoted to RWREs on Galton-Watson trees, our main objective being to study the properties of the limiting velocity and the scaling of such process. Our main focus will be on biased random walks in a wide variety of settings: on Galton-Watson trees with leaves, without leaves and we will also discuss the case of randomly chosen biases.

Then, in Section~\ref{sec_percomain}, we move on to the biased random walks on percolation clusters and present the latest results concerning the long term behavior of such walks (speed, scaling limits ...),

Finally, we shall give an overview of related results on critical structures, namely the infinite incipient cluster (or I.I.C.) on trees or the lattice, in Section~\ref{sec_crit}.

Most of the work presented in these notes was developed in the last five years, hence many questions remain open.  We will include at the end of most sections a list of open problems related to the material presented in the corresponding section.

At the end of the paper, we include a brief appendix  covering some basic facts on sums of heavy-tailed i.i.d.~random variables, a key element for understanding trapping in directionally transient RWREs.

\section{One dimensional models}\label{sect_1d}

Our first goal in this notes is to introduce the key toy model for understanding RWREs in the directionally transient regime. We will show how this elementary model is representative of more complex models, namely one-dimensional RWREs.

\subsection{The totally directed Bouchaud Trap Model}

We introduce a sequence  $(\tau_i)_{i\geq 0}$ of positive i.i.d.~random variables. In the context of trap models, these random variables constitute our environment and represent the strength of the different traps. They are usually assumed to have heavy tails but for the time being we will not make any further assumptions.

The totally directed Bouchaud Trap Model is a nearest-neighbor continuous-time Markov process $(X_t)_{t\geq 0}$ with state space $\Z^+$. Its dynamics are very simple, $X_t$ spends at a site $x$ an exponentially distributed time of mean $\tau_x$ (we recall that $\tau_x$ was chosen randomly) and then jumps to its right neighbor. Our goal is to see how fast $X_t$ moves.

The analysis of this model is pretty elementary. Its interest comes from the fact that several types of biased (or directionally transient) RWREs can be compared to this toy model as will become apparent through these notes.

One key random variable to understand how fast $X_t$ moves is  the time it takes us to make $n$ steps (and thus reach $n$). This quantity can be written
\begin{equation}\label{def_sn}
S(n)=\sum_{i=0}^{n-1} \tau_{i} {\bf e}_i
\end{equation}
 where $({\bf e}_i)_{i\geq 0}$ is a family of i.i.d.~ mean-one exponentially distributed random variables. 
 
It is simple to notice that
\begin{equation}\label{aboveq}
X_t= S^{-1}(t),\qquad \text{for all $t\geq 0$},
\end{equation}
 where the right-continuous inverse of an increasing function $\phi$ is defined by $\phi^{-1}(t):=\inf\{u\geq 0 \mid \phi(u)>t\}$. This equation will allow us to derive scaling limits of $X_t$ from those of $S(n)$ by elementary inversion arguments.

The reason why it is convenient to consider $S(n)$ is that it is simply a sum of i.i.d.~random variables. The topic of sums of i.i.d.~random variables is very well understood and thus it is possible to give a complete picture of all the possible behaviors for scaling limits of $S(n)$. Since a lot of the results in these notes boil down to sums of i.i.d.~random variables, it is very important for the reader to be familiar with the subject. We refer the reader to the Appendix for a brief overview of certain key facts related to sums of i.i.d.~heavy-tailed random variables.

The most trivial case is when ${\bf E}[\tau^2]<\infty$ in which case we have the law of large numbers and a central limit theorem for $S(n)$. Denoting $v:= \frac 1 {{\bf E}[\tau {\bf e}]} $ (for velocity), we have
\[
\lim_{n\to \infty} \frac{S(n)}n = v^{-1} \text{ a.s., and } \lim_{n\to \infty} \frac{ S(n)- nv^{-1}}{ \text{Var}(\tau {\bf e})\sqrt n }\xrightarrow{(d)} \mathcal{N}(0,1),
\]
which, by a standard inversion argument, translates into
\[
\lim_{t\to \infty} \frac{X_t}t =v\text{ a.s., and } \lim_{t\to \infty} \frac{ X_t-tv }{v^{3/2} \text{Var}(\tau {\bf e})\sqrt t }\xrightarrow{(d)} \mathcal{N}(0,1),
\]

One could also easily prove that the processes  
\[
(S_t^{(N)}, 0\leq t\leq T):=(N^{-1/2}(S(\lfloor tN \rfloor)- tN/v), 0\leq t\leq T),
\]
 and 
 \[
  (X_t^{(N)}, 0\leq t\leq T):=(N^{-1/2}(X_{tN}-vtN), 0\leq t\leq T)
  \]
  both  converge as $N$ goes to infinity to  Brownian motions in the Skorohod topology, see Chapter 3 of~\cite{EK}.

In cases where ${\bf E}[\tau^2]=\infty$, we start feeling the effects of traps and witnessing anomalous (non-Gaussian limiting behavior) and we will from now on solely focus on this case. 

To simplify the notations we will assume that ${\bf P}[\tau \geq t]\sim t^{-\alpha}$ for $\alpha \in (0,2]$. Although this might seem like an arbitrary choice, up to a scaling by a constant, this covers all tails of the form $Ct^{-\alpha}$. Most cases appearing in models of biased RWREs fall into that category and, thus, the results in next sections are the ones we expect to witness for many biased RWREs where trapping occurs.

\begin{remark}
Up to some minor changes in the centering and the scaling in the coming results, all of them can be carried over to the case where ${\bf P}[\tau \geq t]\sim L(t)t^{-\alpha}$ where $\alpha \in (0,2]$ and $L(t)$ is a slowly varying function (see the Appendix, definition~\ref{def_sec_varlente}, for a precise definition).
\end{remark}

\subsubsection{The case $\alpha\in (0,1)$ }

In this section, we assume ${\bf P}[\tau \geq t]\sim t^{-\alpha}$ for $\alpha \in (0,1)$. In this case, the trapping is very strong and we actually have ${\bf E}[\tau]=\infty$, this means that even the law of large numbers fails for $S(n)$. Let us state what kind of results can be obtained in this case.

\vspace{0.4cm}

{\it Stable laws as scaling limits}

\vspace{0.4cm}

First, we see that $X_t$ has zero-speed and should be rescaled as $n^{\alpha}$. The proper rescaling leads to stable scaling limits:
\begin{equation}\label{zd_1}
\frac{S(n)}{n^{1/\alpha}} \xrightarrow{(d)} \mathcal{S}_{\alpha}^{ca} \text{ and } \frac{X_t}{t^{\alpha}} \xrightarrow{(d)} (\mathcal{S}_{\alpha}^{ca})^{-\alpha},
\end{equation}
where $\mathcal{S}_{\alpha}^{ca}$ is a completely asymmetric positive stable random variable of index $\alpha$ (see the Appendix, Section~\ref{def_stable_law_section} for a more precise definition).

\vspace{0.4cm}

{\it Stable  subordinators as limiting processes}

\vspace{0.4cm}

 These results carry over to the process level. Indeed, if we consider, for $T>0$ fixed, the processes
\begin{equation}\label{zd_2}
 (S_t^{(N)}, 0\leq t\leq T):=(N^{-1/\alpha}S(\lfloor tN \rfloor), 0\leq t\leq T),
\end{equation}
and
\[
 (X_t^{(N)}, 0\leq t\leq T):=(N^{-\alpha}X_{tN}, 0\leq t\leq T),
 \]
  then the law of $S_t^{(N)}$ defined on the space $D([0,T])$ of c\`adl\`ag functions from $[0,T]$ to $\R$ equipped with the Skorohod $M_1$-topology (see~\cite{whitt} for a detailed account on this topology) converges weakly to the distribution of 
\begin{equation}\label{zd_3}
    ( S_t, 0\leq t\leq T),
\end{equation}
where $S_t$ is an $\alpha$-stable subordinator (see the Appendix, Section~\ref{def_levyproc}, for a precise definition). Moreover, the law of of $X_t^{(n)}$ defined on $D([0,T])$ equipped with the uniform topology converges weakly to the law of 
\begin{equation}\label{zd_4}
  ( Z_t, 0\leq t\leq T),
  \end{equation}
   where $Z_t$ is the inverse of an $\alpha$-stable subordinator.

\vspace{0.4cm}

{\it Aging properties}

\vspace{0.4cm}

Aging is one of the main paradigms in dynamics in random media. It appears in the context of dynamics on spin glasses~\cite{BABC}  (see~\cite{Bouchaud}  for a physical overview of spin glasses) as well as in the random energy model under Glauber dynamics~\cite{bbgtwo}) or in parabolic Anderson model (see~\cite{MOS}). A dynamic satisfies aging if the right time scale to witness a significant change in the system is of the order of the `age' of the system (i.e.~the time during which the dynamics has run). 

In our context, this can be expressed in the following property: for all $a,b>0$ with $a<b$, we have
\begin{equation}\label{zd_5}
\lim_{t\to \infty} \PR[X_{at}=X_{bt}] =\frac{\sin{\alpha \pi}}{\pi} \int_0^{a/b} y^{\alpha-1}(1-y)^{-\alpha}dy,
\end{equation}
where the integral on the right side can be rewritten as $P[\text{ASL}_{\alpha}\in [0,a/b]]$, where $\text{ASL}_{\alpha}$ denotes the generalized arcsine distribution with parameter $\alpha$. The arcsine law is a distribution on $[0,1]$.

\subsubsection{The case $\alpha \in(1,2)$}\label{alpha12}

In this section we assume ${\bf P}[\tau \geq t]\sim t^{-\alpha}$ for $\alpha \in (1,2)$. In this case we still have ${\bf E}[\tau {\bf e}]<\infty$, and thus the law of large numbers holds
\[
\lim_{n\to \infty} \frac{S(n)}n = v^{-1} \text{ a.s., and }\lim_{t\to \infty} \frac{X_t}t =v \text{ a.s.},
\]
where $v:= \frac 1 {{\bf E}[\tau {\bf e}]}$.

The trapping is actually felt in the fluctuations of the walk.

\vspace{0.4cm}

{\it Stable laws in the fluctuations}

\vspace{0.4cm}

\[
\frac{S(n)-nv^{-1}}{n^{1/\alpha}} \xrightarrow{(d)} \mathcal{S}_{\alpha}^{ca} \text{ and } \frac{X_t- vt }{t^{1/\alpha}} \xrightarrow{(d)} -v^{1+1/\alpha}\mathcal{S}_{\alpha}^{ca},
\]
where $\mathcal{S}_{\alpha}^{ca}$ is a completely asymmetric zero-mean stable random variable of index $\alpha$ (see the Appendix, Section~\ref{def_stable_law_section} for a more precise definition).

\vspace{0.4cm}

{\it Stable totally asymmetric L\'evy processes as limiting processes for the fluctuations}

\vspace{0.4cm}

For $T>0$ fixed, if we define the processes
\[
 (S_t^{(N)}, 0\leq t\leq T):=(N^{-1/\alpha}(S(\lfloor tN \rfloor)-tNv^{-1}), 0\leq t\leq T),
\]
and
\[
 (X_t^{(N)}, 0\leq t\leq T):=(N^{-1/\alpha}(X_{tN}-v tN), 0\leq t\leq T),
 \]
  then the law of $S_t^{(N)}$ defined on the space $D([0,T])$ of c\`adl\`ag functions from $[0,T]$ to $\R$ equipped with the Skorohod $M_1$-topology converges weakly to the distribution of 
\[
    ( S_t, 0\leq t\leq T),
\]
where $S_t$ is an $\alpha$-stable totally asymmetric L\'evy process (see the Appendix, Section~\ref{def_levyproc}, for a precise definition). Moreover, the law of $X_t^{(n)}$ defined on $D([0,T])$ equipped with the uniform topology converges weakly to the law of 
\[
  (-\bigl(v^{1+\frac 1{\alpha}}S_t\bigr), 0\leq t\leq T),
\]
   where $S_t$ is an $\alpha$-stable asymmetric L\'evy process (with the same law as the one appearing in the previous equation).

\begin{remark} The cases where $\alpha=1$ and $\alpha=2$ can also be treated (see~\cite{Petrov} or~\cite{GK} for more background on sums of i.i.d.~random variables). We chose not to cover them in these notes, because we will essentially not discuss these cases in the case of RWREs. \end{remark}

\subsection{The biased Bouchaud Trap Model}\label{sect_btm}

Let us define the biased Bouchaud Trap Model. Again we choose a sequence $(\tau_i)_{i\in \Z}$ of positive i.i.d.~random variables, and assume that 
\[
\lim_{x\to \infty} x^{\alpha} {\bf P}[\tau\geq x]=1,
\]
for some $\alpha \in (0,1)$.

For any $\beta>1$, the $\beta$-biased Bouchaud trap model is defined in the following manner: started at $0$, the walk $X_t$ spends at any site $x$ an exponentially distributed time of mean $\tau_x$ and then jumps to the right with probability $\beta/(\beta+1)$ and to the left with probability $1/(\beta+1)$.

\begin{remark}\label{def_snbias}The totally directed Bouchaud trap model is the particular case $\beta=\infty$. \end{remark}

We choose to view $X$ as a time change of a discrete biased random walk on $\Z$. For this, we denote by
\begin{enumerate}
\item $Y_n$  a  $\beta$-biased random walk on $\Z$ recording the successive locations of $X_t$
\item  $S(n)$ the time of the $n$-th jump of $X_t$ which can be written $S(n)=\sum_{i=0}^{n-1} \tau_{Y_i} {\bf e}_i$ where $({\bf e}_i)_{i\geq 0}$ is a family of i.i.d.~ mean-one exponentially distributed random variables.
\end{enumerate}

In this context, $Y_n$ and $S(n)$ are respectively called the embedded process and the clock process.

We emphasize that we changed our definition of $S(n)$ compared to the previous section (we replaced $\tau_i$ by $\tau_{Y_i}$), although in the case $\beta=\infty$ these definitions coincide since $Y_i=i$. It is then simple to notice that
\[
X_t= Y_{S^{-1}(t)},\qquad \text{for all $t\geq 0$},
\]
 where the right-continuous inverse of an increasing function $\phi$ is defined by $\phi^{-1}(t):=\inf\{u\geq 0 \mid \phi(u)>t\}$.

\begin{remark} One of the purposes of this transformation is to notice the link with biased random walks in random environment. The bias is encoded in the embedded process and the random environment is felt in the clock process.  \end{remark}

The key difficulty is the interplay between the past trajectory and the clock process, in particular $S(n)$ is no longer a sum of i.i.d.~random variables. In this context, the time spent in a given trap $i$ is determined by two elements: the number of visits to that trap and the time spent during each of those visits. Although for two distinct traps the time spent during visits in those traps is still independent the number of visits to two adjacent traps is very much correlated. One could then be led to believe that the results obtained through classical theorems on sums of i.i.d.~random variables will fail.

Nevertheless, it turns out that  the following result is true, see~\cite{Zindy}.
\begin{theorem}
Let us consider trapping times $(\tau_i)_{i\geq 0}$ which are positive i.i.d.~random variables such that $\lim_{x\to \infty} x^{\alpha}P[\tau\geq x]=1$, for some $\alpha \in (0,1)$. Then the  $\beta$-directed Bouchaud trap model has stable limit laws, stable  subordinators as limiting processes and verifies the aging property.
\end{theorem}

More precisely, the properties listed at~(\ref{zd_1}), (\ref{zd_2}), (\ref{zd_3}) and~(\ref{zd_4}) remain true in the case of the $\beta$-directed Bouchaud trap model. The only difference is that the limiting $\alpha$-stable laws (or subordinators) have different parameters (see the Appendix, Section~\ref{def_levyproc}, for precise definitions). Furthermore, the property~(\ref{zd_5}) is conserved exactly as it is.

This result raises a natural question: why does the model not feel the dependencies?

\begin{remark} Similar results to Section~\ref{alpha12}, should be true for this model when $\alpha \in (1,2)$. \end{remark}
                                                                                                                                                                                                                                                                                                                                                                                                                                                                                                                                                                                                                                                                                                                                                                                                                                                                                                                                                                                                                                                                                                                                                                                                                                                                                                                                                                                                                                                                                                                                                                                                                                                                                                                                                                                                                                                                                                                                                                                                                                                                                                                                                                                                                                                                                                                                                                                                                                                                                                                                                                                                                                                                                                                                                                                                                                                                                                                                                                                                                                                                                                                                                                                                                                                                                                                                                                                                                                                                                                                                                                                                                                                                                                                                                                                                                                                                                                                                                                                                                                                                                                                                                                                                                                                                                                                                                                                                                                                                                                                                                                                                                                                                                                                                                                                                                                                                                                                
\subsubsection{Why do trap models ignore many  dependencies?}\label{trap_dep}

Let us consider trapping times $(\tau_i)_{i\geq 0}$ which are positive i.i.d.~random variables such that 
\[
\lim_{t\to \infty} t^{\alpha}{\bf P}[\tau\geq t]=1,
\]
 for some $\alpha \in (0,1)$.

Our goal is to intuitively explain why a directed trap model with such trapping times is insensitive to many correlations: why is it that $S(n)$ (defined after Remark~(\ref{def_snbias})) still behaves as a sum of i.i.d.~random variables?

To compute the hitting time $S(n)$, the interval $[0,n]$ is essentially the only one relevant for us. Indeed the directed nature of the walk makes us unlikely to visit many sites at the left of the origin (at most $C\ln n$). An elementary computation shows that the deepest traps  (i.e.~those of maximal mean waiting time) encountered in the interval $[0,n]$ is of the order $n^{1/\alpha}$. Fix $\epsilon>0$, the previous brings us to the following observations.

\begin{itemize}
\item Before $S(n)$, we are only visiting traps of depth (i.e.~mean waiting time) less than $n^{(1 +\epsilon)/\alpha}$.
\end{itemize}

Also, we can actually estimate the number of deep traps encountered.

\begin{itemize}
\item Before $S(n)$, we visit roughly $n^{\epsilon}$ traps of depth larger than $n^{(1-\epsilon)/\alpha}$. They are also spatially very well separated (at typical distance $n^{1-\epsilon}$). Since the walk is unlikely to backtrack (go towards the left) for more than $C\ln n$ steps, the behavior of the walk (number of visits ...) in two different large traps is essentially independent.
\end{itemize}

Finally, small traps are irrelevant as far as limiting behaviors are concerned.

\begin{itemize}
\item Before $S(n)$, the total time spent in all traps of depth less than $n^{(1-\epsilon)/\alpha}$ is $o(n^{1/\alpha})$.
\end{itemize}

Hence
\[
S(n) \approx \sum_{i=0}^{n^{\epsilon}} T_n^{(i)},
\]
where the $T_n^{(i)}$ are i.i.d.~and represent the total time spent in traps conditioned to be larger than $n^{(1-\epsilon)/\alpha}$. This setting is not the standard one of sums of i.i.d.~random variables. Indeed, the law of the terms is $n$-dependent. Nevertheless this object is known as a triangular array and the limiting behavior of triangle arrays is very well understood (see~\cite{Petrov}). It is by using limiting theorems on such triangular arrays that we can manage to obtain stable scaling limits. 

\subsection{The one dimensional random walk in random environment}

Up to this point we have only considered toy models, in which traps were introduced ad hoc, by forcing waiting times at different sites. This allowed us to understand which limiting results we can expect when trapping is involved and at the same time it gave us an opportunity to describe the way trap models work.

We will present a random dynamic in random media which does not have the trapping encoded in its definition but yet exhibits the same type of limiting behaviors: one dimensional random walks in random environments.

One-dimensional random walks in random environment to the nearest
neighbors were introduced in the sixties in order to give a
model of DNA replication. Mathematically, these models were initially studied in the seventies and are still yielding new results up to the present day.

 \subsubsection{The model}
 
Let $\omega:=(\omega_i, \, i \in \Z)$ be a family of i.i.d. random
variables taking values in $(0,1)$ defined on $\Omega,$ which stands
for the random environment. Denote by ${\bf P}$ the distribution of
$\omega$ and by ${\bf E}$ the corresponding expectation. Conditioning on
$\omega$ (i.e. choosing an environment), we define the random walk
in random environment $(X_n, \, n \ge 0)$ as a nearest-neighbor
random walk on $\Z$ with transition probabilities given by $\omega$:
$(X_n, \, n \ge 0)$ is the Markov chain satisfying $X_0=0$ and for
$n \ge 0,$
\[
P^\omega \left[ X_{n+1} = x+1 \, \mid \, X_n =x\right] = \omega_x = 1-
P^\omega \left[ X_{n+1} = x-1 \, \mid \, X_n =x\right].
\]

\begin{figure}
\centering 
\epsfig{file=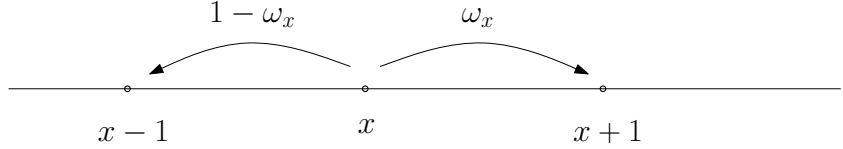, scale=0.7}
\caption{Transition probabilities for the one-dimensional random walk in random environment.}
\end{figure}

We denote by $P^{\omega}$ the law of $(X_n, \, n \ge 0)$
and $E^{\omega}$ the corresponding expectation. It is usually called quenched law. Also we denote by $\PR$
the joint law of $(\omega,(X_n)_{n \geq 0})$, which is referred to as averaged or annealed law. 

 \subsubsection{The results}

In the study of one-dimensional random walks in random environment,
an important role is played by the sequence of variables 
\[
\rho_i:= \frac{1-\omega_i}{\omega_i}, \qquad i \in \Z,
\]
since their law will allow to see when the walk is transient, recurrent, has positive speed, or not.

We assume that ${\bf E}[\ln \rho_0]$ is well defined (possibly infinite). Then we obtain a criterion for transience and recurrence, see~\cite{Solomon}.
\begin{theorem}
\begin{enumerate}
\item If ${\bf E}[\ln \rho_0]<0$ (respectively >0) then the walk is transient and 
\[
\lim_{n\to \infty} X_n = \infty\ (\text{resp. $- \infty$}) \qquad \PR\text{-a.s.}.
\]
\item If ${\bf E}[\ln \rho_0]=0$, then the walk is recurrent and  
\[
\limsup_{n\to \infty} X_n = \infty \text{ and } \liminf_{n\to \infty} X_n =- \infty \qquad \PR\text{-a.s.}.
\]
\end{enumerate}
\end{theorem}

An explanation for the result of this theorem is given at the beginning of Section~\ref{sect_1dbtm}.

Our goal being to focus on the transient case, so we shall henceforth assume  that ${\bf E}[\ln \rho_0]<0$ (which is possible by symmetry). If the reader wants more information on the recurrent case, he/she should consult~\cite{Sinai}. Since we assume directional transience to $+\infty$ we can introduce, for $n\geq 0$, 
\begin{equation}\label{zdefdelta}
\Delta_n=\inf\{i\geq 0, X_i=n\}<\infty,
\end{equation}
 the hitting time of $n$.

We may decide whether the walk has positive speed or zero velocity, see~\cite{Solomon}.
\begin{theorem}[Solomon-1975] 
If we have ${\bf E}[\ln \rho_0]<0$, then 
\[
\lim_{n\to \infty} \frac{\Delta_n}n = v^{-1}\text{ and } \lim_{n\to \infty} \frac{X_n}n = v \qquad \PR\text{-a.s.,}
\]
with
\[
v=\begin{cases} \frac{1-{\bf E}[\rho_0]}{1+{\bf E}[\rho_0]}>0 & \text{if } {\bf E}[\rho_0]<1 \\
                       0 & \text{ otherwise.}
                       \end{cases}
\]
\end{theorem}

The key fact here is that it is possible to obtain directional transience towards $+\infty$ and yet have zero speed which is indicative of trapping. This means we should be expecting stable laws in the limit, and indeed
\begin{theorem}\label{KKS}
Let us assume that 
\begin{enumerate}
 \item $-\infty \leq {\bf E}[\ln \rho_0]<0$,
 \item there is a $0<\alpha<2$ such that ${\bf E} \left[\rho_0^{\alpha}\right] = 1$ 
  and ${\bf E} \left[  \rho_0^{\alpha} \ln^+\rho_0  \right]<\infty,$
 \item  the distribution of $\ln \rho_0$ is non-lattice (that is, the $\Z$-linear span of $\text{supp}(\ln \rho_0)$ is dense in $\R$),
\end{enumerate}
then 
\begin{itemize}
\item if $1<\alpha <2$, letting $v:= \frac{1-{\bf E}[\rho_0]}{1+{\bf E}[\rho_0]} $, we have
\[
\frac{\Delta_n-nv^{-1}}{n^{1/\alpha}} 
\stackrel{law}{\longrightarrow}\, \mathcal{S}_{\alpha}^{ca}  \text{ and }
\frac{X_n-nv}{n^{1/\alpha}} 
\stackrel{law}{\longrightarrow}\, -v^{1+1/\alpha}\mathcal{S}_{\alpha}^{ca},
\]
\item if $0<\alpha <1$, then
\[
\frac{\Delta_n}{n^{1/\alpha}} 
\stackrel{law}{\longrightarrow}\, \mathcal{S}_{\alpha}^{ca} \text{ and }
\frac{X_n}{n^{\alpha}}
\stackrel{law}{\longrightarrow}\,   (\mathcal{S}_{\alpha}^{ca})^{-\alpha},
\]
\end{itemize}
where  $\mathcal{S}_{\alpha}^{ca} $ is a completely asymmetric stable law of index $\alpha$.
\end{theorem}

A more complete result, covering the case $\alpha=1$ and $\alpha=2$, can be found in~\cite{KKS}. The original proof of Kesten, Kozlov and Spitzer gives a very complete result using a short proof related to branching processes. Although efficient this argument is very specific to $\Z$ and could not be carried over to different underlying graphs.

An alternate proof of the previous result was obtained in~\cite{ESZ} and~\cite{ESZ1} (for $\alpha\in (0,1)$) and later completed in~\cite{ESTZ} (for $\alpha \in [1,2)$). This alternate proof used the similarities of one-dimensional RWREs and the biased Bouchaud trap model. To the best of our knowledge, it was the first successful attempt to link directly RWREs to trap models. The works~\cite{ESZ} and~\cite{ESZ1} actually predate and inspired the work on biased trap models~\cite{Zindy} that we presented earlier.

\begin{remark}  The methods of~\cite{ESZ} and~\cite{ESZ1} are sufficient to prove the convergence in Theorem~\ref{KKS} at the process level. \end{remark}

\begin{remark} \label{extra_info} Under the hypotheses of Theorem~\ref{KKS}, it is known, see~\cite{P} and~\cite{PetZ}, that quenched scaling does not exist for the one-dimensional RWRE when $\alpha\in (0,2)$. 

However, as a by-product of these methods in~\cite{ESZ} and~\cite{ESZ1}, it is possible to obtain a description of the long time behavior of the quenched law of the hitting time associated to the walk (see~\cite{ESZ3} and~\cite{ESTZ1}). We mention this result was proved in parallel in independent works of~\cite{PS} and~\cite{DG}.  \end{remark}

A key benefit of this trap model approach is that it allows one to obtain the following aging result when $\alpha\in (0,1)$, which was previously unknown, see~\cite{ESZ3}.
\begin{theorem}[Enriquez, Sabot, Zindy - 2009]\label{aging_esz}
Let us assume that 
\begin{enumerate}
 \item $-\infty \leq {\bf E}[\ln \rho_0]<0$,
 \item there is a $0<\alpha<1$ such that ${\bf E} \left[\rho_0^{\alpha}\right] = 1$ 
  and ${\bf E} \left[  \rho_0^{\alpha} \ln^+\rho_0  \right]<\infty,$
 \item  the distribution of $\ln \rho_0$ is non-lattice,
\end{enumerate}
then, for any $\eta>0$ and any $h>1$, we have
 \[
 \lim_{t\to \infty} \PR[\abs{X_{th}-X_t} \leq \eta \ln t]=\frac{\sin{\alpha \pi}}{\pi} \int_0^{1/h} y^{\alpha-1}(1-y)^{-\alpha}dy.
 \]
 \end{theorem}

 An other, even more important, consequence of this work is that it provided methods that could later be used to analyze RWREs on different graphs such as trees and $\Z^d$ ($d\geq 2$).

 \subsubsection{Parallel between one dimensional RWREs and the biased BTM}\label{sect_1dbtm}
 
The parallel one dimensional RWREs and the biased BTM can be made by considering an object $V(x)$ called the potential that was introduce in 1982 by Sinai~\cite{Sinai}.
\begin{equation}\label{potential}
V(x) :=\left\{\begin{array}{ll} 
          \sum_{i=1}^x \ln \rho_i, & \text{if } x \geq 1, \\
 0\vphantom{\sum^N}, &  \text{if }  x=0, \\
-\sum_{i=x+1}^0 \ln \rho_i, &\text{if } x\leq -1,
\end{array}\right.
\end{equation}

This potential gives a nice intuitive vision of the transition probabilities. Loosely speaking it could be seen as an altitude profile on which the walker will move, favoring directions going downhill. Mathematically, we see that 
\begin{itemize}
\item $V(x)$ is a random walk with i.i.d.~steps of law $\ln \rho_0$. This means it can be analyzed very precisely.
\item $V(x)$ is linked to an invariant measure $\pi$ of the walk in the environment $\omega$ by setting $\pi(x)=e^{-V(x)}+e^{-V(x-1)}$. This explains why the walk is attracted to parts where $V(\cdot)$ is small.
\end{itemize}

\begin{figure}
\centering 
\epsfig{file=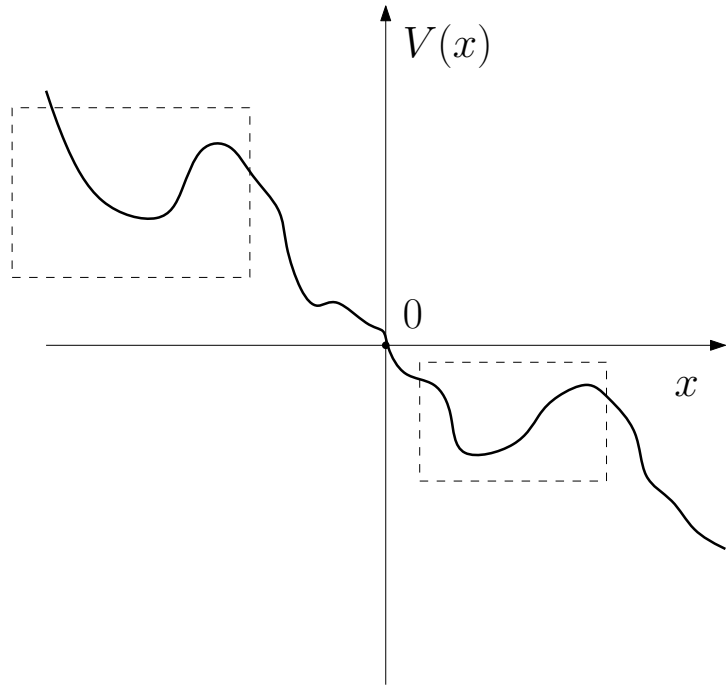, scale=0.8}
\caption{The typical potential when the walk is transient towards $+\infty$.}
\end{figure}

Let us  detail the strategy used in~\cite{ESZ} and~\cite{ESZ1} to link one-dimensional RWREs to the biased BTM using the potential $V$. In Figure~6, we have singled out parts of the environment where the walk needs to walk uphill against the potential to reach $+\infty$. These valleys create obstacles that prevent the walk from going easily to $+\infty$. These obstacles have a random height $H$ which can be shown to verify
\begin{equation}\label{1d_height}
\PR[H \geq n] \sim C \exp(- \alpha n ),
\end{equation}
where $\alpha$ appears in~Theorem~\ref{KKS} by a result of~\cite{Iglehart}.

When in a valley of height $H$, the walk will need a certain time $T_{\text{exit}}$ to exit it. To simplify the argument, we use the approximation $T_{\text{exit}}=e^H{\bf e}$ where ${\bf e}$ is an exponential random variable of parameter 1. A reader interested in deeper details on this issue might want to look at Section~\ref{sect_approx} and Remark~\ref{1d_trapex} for a deeper discussion on a similar problem in a different context. Using the tail estimate on $H$ we obtain
\[
\PR[T_{\text{exit}}\geq n]\sim C n^{-\alpha},
\]
which is the heavy-tailed random variable we were expecting to see when $\alpha\in (0,1)$.

The parallel with the biased BTM should now be pretty clear. A few issues might still bother the reader:
\begin{enumerate}
\item The embedded walk is no longer a biased random walk. While this is true, we only really need for the embedded walk to go against the direction of transience for small distances  (at most $C\ln n$ before time $\Delta_n$, see~(\ref{zdefdelta}) for the definition). This is still verified in the model.
\item The trap depth are no longer independent, indeed many neighboring sites are in the same valley and share similar exit times. This issue can be resolved by considering the entire valleys as one trap (or a site in the biased BTM) since valleys are typically small (at most $C\ln n$ before  time $\Delta_n$) they appear as points when rescaled in an environment of length $n$.  One can then argue that the heights and geometry of the largest traps is roughly independent since they are located far from each other.
\end{enumerate}

\section{Biased random walk on supercritical trees}\label{sect_arbre_model}

\subsection{The model}

The model of biased random walk on supercritical Galton-Watson trees was initially studied in~\cite{LPP}. As discussed in the introduction, it has an interesting phenomenology and is a natural setting for discussing trapping issues.

\subsubsection{Definition}\label{def_model}

Consider a supercritical Galton-Watson branching process with generating function ${\bf f}(z)=\sum_{k \geq 0} p_k z^k$, i.e.~the offspring of all individuals are i.i.d.~copies of $Z$, where $ P[Z=k]=p_k$. We assume that the tree is supercritical, i.e. ${\bf m}:= E[Z]={\bf f}'(1)\in (1,\infty)$. If $p_0=0$, the tree has no leaves and otherwise it does.
We denote by $q$ the extinction probability ($q<1$ since the tree is supercritical), which is characterized by ${\bf f}(q)=q$. We know that $q=0$ if, and only if, $p_0=0$. Starting from a single progenitor called root and denoted by $0$, this process yields a random tree $T(\omega)$. If $p_0>0$, we will always condition on the event of non-extinction, so that $T$ is an infinite random tree. We denote by $(\Omega,{\bf P})$ the associated probability space:  ${\bf P}$ is the law of the original tree, conditioned on non-extinction. 

Given such an infinite Galton-Watson tree $T(\omega)$, we consider the
$\beta$-biased random walk started at $0$ as defined in Section~\ref{sect_intro_biased}. For a fixed environment $\omega$, the associated law of the $\beta$-biased random walk is called quenched and is denoted $P^{\omega}$. We define the averaged (sometimes called annealed) law as the semi-direct product $\PR={\bf P} \times P^{\omega}$.

\begin{remark} In a fixed environment $\omega$, this is a reversible Markov chain, with an invariant measure given by $\pi(x)=\sum_{y\sim x} c(x,y)$ where every edge $[\overleftarrow{x}, x] $ of the tree is given a `conductance' $c(\overleftarrow{x},x):=c(x,\overleftarrow{x}):=\beta^{\abs{x}-1}$ (see~\cite{LP} or~\cite{DoyleSnell} for background on electrical networks). 
\end{remark}

For a vertex $u \in T$, we denote by $\abs{u}=d(0,u)$ the distance of $u$ to the root.

\subsubsection{Basic results}

Expanding on the results of Section~\ref{sect_intro_biased}, we have, by~\cite{lycap}, the following  recurrence and transience criterion.
\begin{theorem}\label{lycap}
Let $X_n$ be a $\beta$-biased random walk on a Galton-Watson tree with law ${\bf P}$.  We have the following
 \begin{itemize}
\item if $\beta<1/ {\bf m}$, then $X_n$ is strongly recurrent $P^{\omega}$-a.s.~ for ${\bf P }$-almost all $\omega$.
\item if $\beta = 1/{\bf m}$, then $X_n$ is null recurrent $P^{\omega}$-a.s.~ for ${\bf P }$-almost all $\omega$
\item if $\beta > 1/{\bf m}$, then $X_n$ is transient $P^{\omega}$-a.s.~ for ${\bf P }$-almost  all $\omega$ which means that $\lim_{n\to \infty}  \abs{X_n} =\infty$,
\end{itemize}
\end{theorem}

In the sequel we will, unless stated otherwise, assume that $\beta >1/{\bf m}$ which means that we have transience of $\abs{X_n}$ towards $+\infty$.

The next natural question, that of the behavior of the speed, was investigated in~\cite{LPP}.
\begin{theorem}\label{gw_lppspeed}
We have
\[
\lim_{n\to \infty} \frac{\abs{X_n}}n = v(\beta,{\bf P}), \qquad \text{$P^{\omega}$-a.s.~ for ${\bf P }$-a.s. all $\omega$}
\]

Furthermore 
\begin{enumerate}
\item if $p_0=0$, then $v>0$,
\item if $p_0>0$, then $v>0$ if $\beta\in (1/{\bf m}, \beta_c)$ and $v=0$ if $\beta \geq \beta_c$ where $\beta_c=1/{\bf f}'(q)$.
\end{enumerate}
\end{theorem}

\begin{remark} We emphasize that $\beta_c >1$ (see Remark~\ref{fprimeq}). This means that the zero-speed regime only occurs if $\beta>1$, that is if the walk is outwardly-biased.\end{remark}

This means that we can answer  the second question asked at the end of the introduction (just before Section~\ref{sect_plan}) affirmatively: the speed has a sharp phase transition from positive speed to zero speed.

 As previously mentioned this result is due to the existence of dead ends in the tree creating trapping structure for a strongly biased walk, see Figure 1. In the next section, we will introduce the Harris-decomposition for Galton-Watson trees. This construction  allows us to describe precisely the structure of traps. It is key to understand the link between biased random walks on Galton-Watson trees and traps models. This subject will be treated in more details in Section~\ref{sect_arbre_trap}.

\subsubsection{The Harris-decomposition of Galton-Watson trees}\label{sect_arbre_harris}

 Set 
\begin{equation}
\label{notationgh}
{\bf g}(s)=\displaystyle{\frac{{\bf f}((1-q)s+q )-q }{1-q}} \text{ and } {\bf h}(s)=\frac{{\bf f}(qs)}q.
\end{equation}

\begin{figure}
\centering 
\epsfig{file=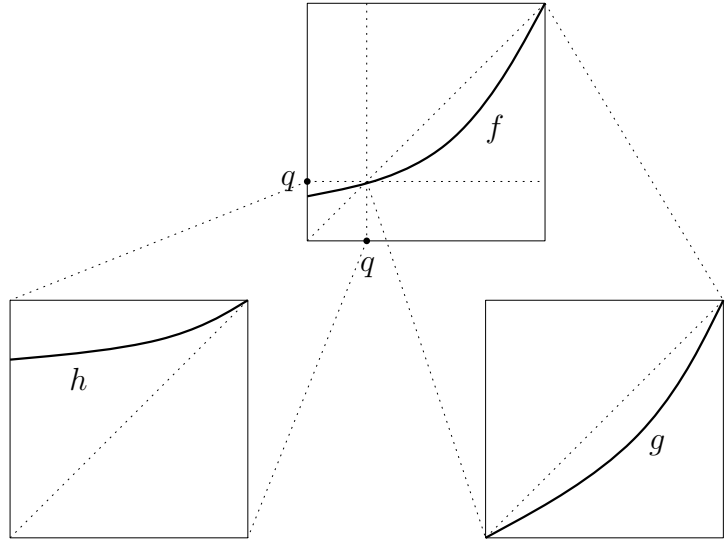, scale=0.7}
\caption{The generating functions in the Harris-decomposition}
\end{figure}

\begin{remark} \label{fprimeq}The functions ${\bf g}$ and ${\bf h }$ have a nice graphical representation, see Figure~7. Here, we can easily see that ${\bf f}'(q)<1$.\end{remark}

It is known (see~\cite{lycap}), that a ${\bf f}$-Galton-Watson tree (with $p_0>0$) can be generated by
\begin{enumerate}
\item growing a ${\bf g}$-Galton-Watson tree $T_{{\bf g}}$ called the  backbone, where all vertices have an infinite line of descent,
\item attaching on each vertex $x$ of $T_{{\bf g}}$ a random number $N_x$ of subcritical (hence finite) ${\bf h}$-Galton-Watson trees, acting as traps for the biased random walk,
\end{enumerate}
where $N_x$ has a distribution depending only on $\text{deg}_{T_{{\bf g}}}(x)$ and given $T_{{\bf g}}$ and $N_x$ the traps are i.i.d. (see \cite{lycap} for details).

\begin{figure}
\centering 
\epsfig{file=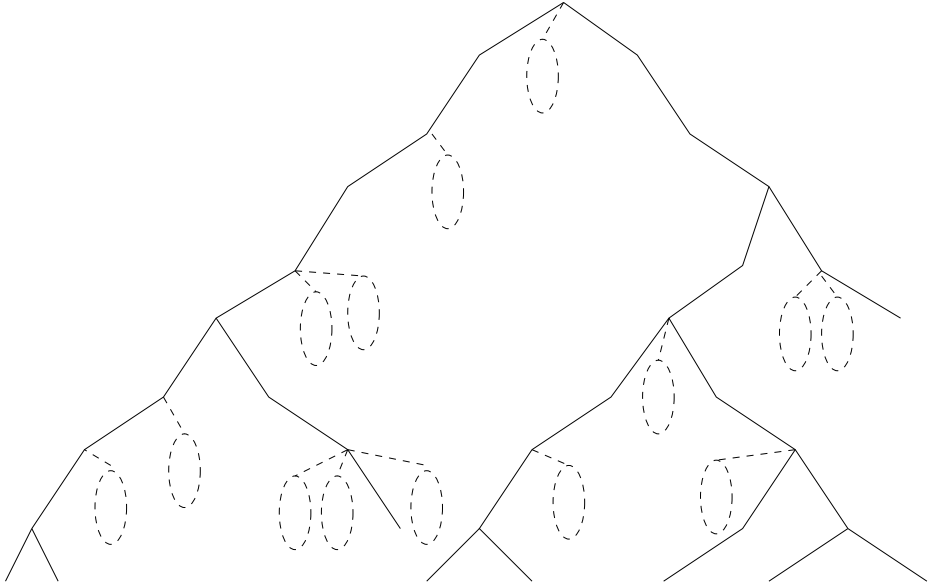, scale=0.7}
\caption{The Galton-Watson tree is decomposed into the backbone (solid lines) and the traps (dashed lines).}
\end{figure}

The Harris-decomposition, see Figure~8, allows us to view the biased random walk on the Galton-Watson tree with leaves in terms very similar to those of the trap models introduced so far. Indeed, one may view the walk on the backbone as the {\bf embedded process} on a leafless ${\bf g}$-Galton-Watson tree with waiting times in  traps given by excursion times in ${\bf h}$-Galton-Watson trees (recall Section~\ref{sect_btm} for the vocabulary associated to trap models).

\subsection{The asymptotic speed}

Let us  turn to the first question at the end of the introduction (just before Section~\ref{sect_plan}), and try to see if, in the case with leaves, it is true that the speed is unimodal, meaning it is first increasing and then decreasing.

At first sight, the embedded process and the waiting times are in competition to determine the value of the velocity.
\begin{enumerate}
\item On the backbone, there are no hard traps for the walk. This would lead one to believe that increasing the bias would increase the speed of the embedded process.
\item On the other hand, increasing the bias increases the time spent in traps and so delays the walk.
\end{enumerate}

The second point can be shown to be true, at least in expectation, using the mean return time formula (a particular case of the commute time formula, see~\cite{commute}). The first point, however, is not as simple of a question as it may seem. It leads us to the following question which was asked in 1996 (see~\cite{LPP} and~\cite{LPP1}).
\begin{question}\label{conj_speed}
Is the speed of a biased random walk on a Galton-Watson tree without leaves increasing (or even non-decreasing) ?
\end{question}

\subsubsection{The asymptotic speed on  trees without leaves}

 At first sight it seems that the speed of a biased random walk (as defined in the introduction) on any leafless tree $T$ should be non-decreasing. Nevertheless this is wrong as some examples from~\cite{LPP1} show. This section is devoted to presenting some of the examples of~\cite{LPP1} and  their behavior. 

\vspace{0.4cm}

{\it Example 1: Binary tree with pipes}

\vspace{0.4cm}

Let $T$ be a binary tree to which we add a unary tree, which we refer to as a pipe, at every vertex, see Figure~9.

\begin{figure}
\centering 
\epsfig{file=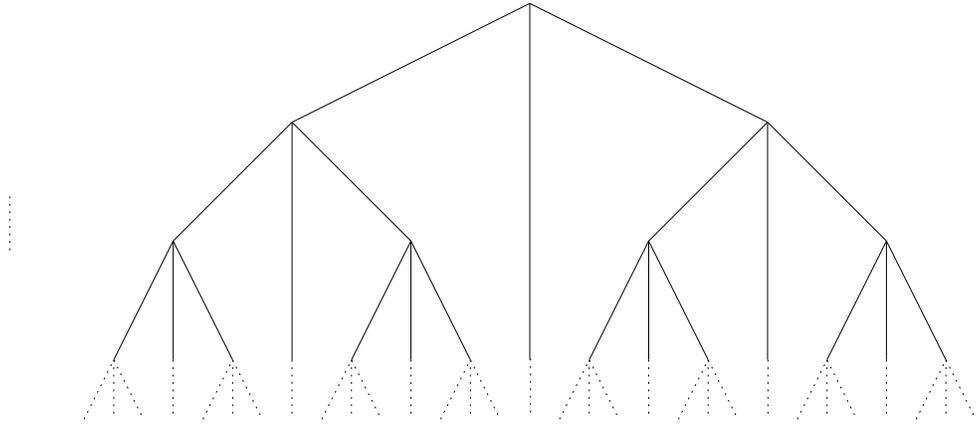, scale=0.7}
\caption{The binary tree with pipes}
\end{figure}

 It is elementary to see that the walk becomes transient for $\beta>1/2$. We can then observe the following.
\begin{enumerate}
\item For $\beta=1/2+\epsilon$, the walk has barely turned transient and its speed should be small: $v(1/2+\epsilon)=O(\epsilon)$.
\item For $\beta=1$, we are looking at the simple random walk which is known to have a (finite) return time with infinite expectation in each pipe. The walk can only go to infinity through the binary part of the tree. However, on its way it will regularly enter pipes and lose a substantial amount of time in them (each excursion having infinite mean). This will mean that $v(1)=0$.
\item For $\beta =3/4$, the walk on the binary part has speed $1/5$ (by comparing it to a $3/2$-biased random walk on $\Z$) and the excursions on the pipes will be brief. This means that $v$ is bounded away from $0$.
\end{enumerate}

For more details on this, we refer the reader to Example 2.1 in~\cite{LPP1}. 

Even though this example clearly shows a mechanism for slowing down the walk whilst increasing the bias, it is slightly dissatisfying. Indeed the cause for the slowdown observed in the second point is that pipes act, de facto, as traps. This makes this example somewhat artificial. In order to prevent such examples, one could impose for the bias to be larger than $1$. This would prevent us from creating hard traps without using leaves. 

\begin{remark}
We emphasize that in the previous example the speed coincides with that of a $\beta$-biased random walk on $\Z$ for $\beta>1$.
\end{remark}

{\it Example 2: The filtering method}

\vspace{0.4cm}

We shall explain how to proceed to create a tree on which the $\beta_1$-biased random walk has a higher speed than the $\beta_2$-biased random walk for $\beta_1<\beta_2$. 

The key fact we use is the following: there exists a tree $T_a$ such that given $\epsilon>0$ there exists  $M$ sufficiently large and two complementary subsets $B_1$ and $B_2$ of the $M$-th level of $T_a$ such that the first visit of the $\beta_i$-biased random walk to level $M$ on $T_a$ first reaches $B_i$ with probability at least $1-\epsilon$. 

This property is verified by almost every tree produced from a Galton-Watson process with mean larger than $\beta_2^{-1}$. We take this fact as a given for our upcoming construction. This property of Galton-Watson trees follows from the fact that on such a tree the harmonic measure of the $\beta_i$-biased random walk are singular, a fact proved in Lemma 5.2 of~\cite{LPP}. The reader may find more details on this in~\cite{LPP} or~\cite{LP}.

Then we apply the following procedure: keep the first $M$ levels of $T_a$, add $n$ (with $n>>M$) levels of a $K_i$-ary tree to each vertex of $B_i$ and we call the resulting tree $T_b$. We choose $K_i$ large enough to ensure that the $\beta_i$-biased random walks are unlikely to revisit the root of a $K_i$-ary tree. Also we choose $\beta_2 K_2>\beta_1 K_1$, so that the  $\beta_2$-biased random walk on the $K_2$-ary tree is faster than the $\beta_1$-biased random walk on a $K_1$-ary tree. 

\begin{figure}
\centering 
\epsfig{file=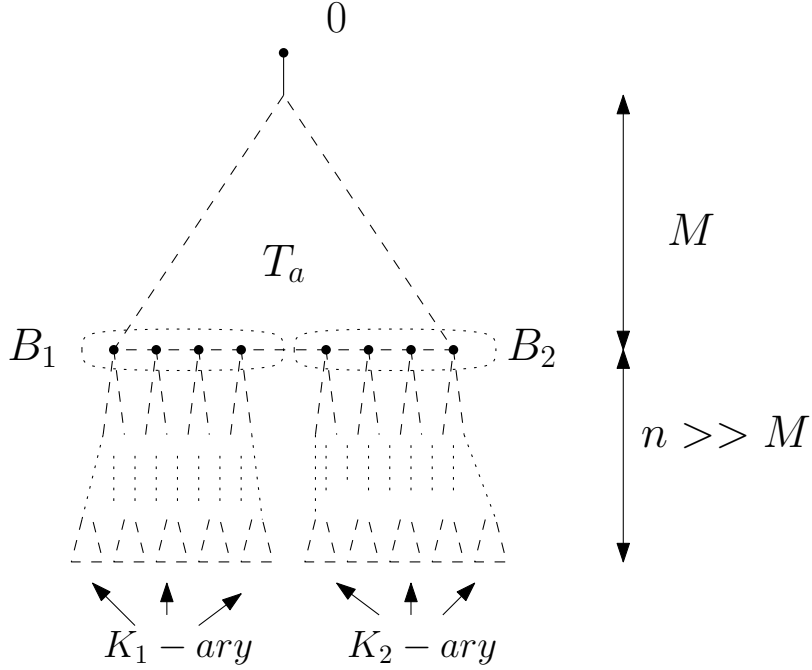, scale=1}
\caption{The construction of the tree $T_b$. Here $K_1$ is much larger than $K_2$.}
\end{figure}

The dynamics on this tree, see Figure~10, are as follows
\begin{enumerate}
\item the $\beta_i$-biased random walk will, w.h.p., be in $B_i$ when it reaches the $M$-th level of $T_b$,
\item from there, the walk will move for a long time on a $K_i$-ary tree, hence, its overall speed will progressively become close to that of the $\beta_i$-biased random walk on a $K_i$-tree.
\end{enumerate}

Essentially, the time spent in the first section of the tree can be made small in comparison to the time spent in the second section, so the total speed is determined mainly by the second part of the journey. In particular the $\beta_2$-biased random walk is faster than the $\beta_1$-biased random walk on $T_b$.

To finish the argument, we take $T_b$, add to each of its leaves a copy of $T_b$ and repeat the procedure constructing an infinite tree $T$. It is clear that on $T$, increasing the bias from $\beta_1$ to $\beta_2$ will slow the walk down.

\begin{remark} The construction of the filtering part of the tree which separates walks of different bias ($T_a$ in the previous example) can also be done in a deterministic manner. Hubert Lacoin has exhibited explicit examples on which the speed can be computed exactly and, indeed, turns out not to be increasing. (Personal communication from Hubert Lacoin).
\end{remark}

Of course none of these trees are Galton-Watson trees (although they are multi-type Galton-Watson trees) so we do not have any counter-examples to Question~\ref{conj_speed}. We merely wanted the reader to appreciate the complexity of the question.

\subsubsection{The asymptotic speed on Galton-Watson trees without leaves}

Recently, some progress has been made in understanding the speed of biased random walks on Galton-Watson trees without leaves. We shall review the main results.

First, one can partially answer the Question~\ref{conj_speed}. Indeed, we have the following result from~\cite{BAFS}
\begin{theorem}\label{BAFS}
The speed $v(\beta)$ of a $\beta$-biased random walk on a Galton-Watson tree without leaves is increasing for $\beta \geq 1160$.

Furthermore, for a $\beta$-biased random walk on a Galton-Watson tree with minimal degree $d:=\min\{k\geq 1, {\bf P}[Z=k] >0\}$, the speed $v(\beta)$ of this biased random walk is increasing for $\beta \geq 1160/d$.
\end{theorem}
 
 \begin{remark} The threshold 1160 is completely arbitrary and could be lowered using the methods of the paper. Nevertheless, the threshold cannot be lowered to any value that would have any particular significance such as $\beta=1$ or $\beta=1/{\bf m}$.
 \end{remark}
 
 The proof  relies on a coupling of three distinct walks: two biased walks ($\beta$ and $\beta+\epsilon$) on Galton-Watson trees and a $\beta$-biased random walk on $\Z$. Once this coupling is in place the proof turns out to be surprisingly elementary.
 
  In some very vague sense, the high bias ensures that the biased walks interact very little with their environment. With the environment out of the picture, the decisive element for a faster speed is the bias. The third walk in the coupling provides an upper bound on the interaction of the walks with the environment.
  
  \begin{remark}\label{rem_multi} The techniques in the proof of Theorem~\ref{BAFS} are somewhat robust. This would lead us to believe that they could be adapted to different problems. In particular, we would think that it should be possible to use these methods to prove that the speed of a $\beta$-biased random walk on a multi-type Galton-Watson tree without leaves that has at most $k$ types is increasing for $\beta \geq \beta_c(k)$ where $\beta_c(k)<\infty$ .
\end{remark}

The coupling  similar to the one introduced in~\cite{BAFS} to prove Theorem~\ref{BAFS} was used in~\cite{NYU} to prove a result on the monotonicity of the speed of biased random walks with respect to the offspring distribution.
\begin{theorem}\label{NYU}
Consider biased random walks on two Galton-Watson trees without 	leaves	having	offspring	distributions	${\bf P}_1$	and	${\bf P}_2$ where ${\bf P}_1$ dominates ${\bf P}_2$ stochastically. Then we have that for a bias that is larger than an explicit threshold depending on ${\bf P}_1$ and ${\bf P}_2$,
\[
v(\beta,{\bf P}_1)\geq v(\beta,{\bf P}_2),
\]
where we used the notations of Theorem~\ref{gw_lppspeed}.
\end{theorem}

\begin{remark}\label{FL} It is natural to expect that using methods of Theorem~\ref{BAFS}, the result of Theorem~\ref{NYU} could be improved  to show that there exists a critical bias $\beta_c$ (independent of ${\bf P_1}$ and ${\bf P_2}$)  above which $v(\beta,{\bf P}_1)\geq v(\beta,{\bf P}_2)$ where ${\bf P_1}$ and ${\bf P_2}$ are supported on positive integers and ${\bf P}_1$ dominates ${\bf P}_2$ stochastically (personal communication by Hubert Lacoin).
\end{remark}

Before moving on to the case of Galton-Watson trees with leaves, we would also like to present a consequence of~\cite{Aidekon} which improves Theorem~\ref{BAFS}  (personal communication from Elie A\"id\'ekon). The result from which the following theorem is derived will be presented in the beginning of the next section.
\begin{theorem}\label{aid1}
The speed $v(\beta)$ of a $\beta$-biased random walk on a Galton-Watson tree without leaves is increasing for $\beta \geq 2$.
\end{theorem}

\begin{remark} Another interesting result about the asymptotic speed on Galton-Watson trees without leaves  is the Einstein relation. This topic will be discussed in Section~\ref{sect_arbre_einstein}. \end{remark}

\subsubsection{The asymptotic speed on Galton-Watson trees with leaves}

A major progress made towards understanding the speed of biased random walks on Galton-Watson trees is due to Elie A\"id\'ekon~\cite{Aidekon}. To state the main result, let us first introduce the random variable $\mathcal{E}_{\infty}:=P_0^{\omega}[X_n\neq 0, \text{ for $n\geq 1$}]$. This random variable is chosen under the (unconditioned) Galton-Watson measure $P$.

In~\cite{Aidekon}, A\"id\'ekon obtains an explicit expression for the speed (we stress that this result is valid on any supercritical Galton-Watson tree, even with leaves).
\begin{theorem}\label{aidekon}
For a biased random walk on a Galton-Watson tree in the positive speed regime, the speed is given by
\[
v(\beta)=\frac{E\Bigl[\frac{(Z-\beta^{-1})\mathcal{E}_{\infty}^{(0)}}{\beta^{-1}-1 +\sum_{i=0}^Z \mathcal{E}_{\infty}^{(i)}}\Bigr]}{E\Bigl[\frac{(Z+\beta^{-1})\mathcal{E}_{\infty}^{(0)}}{\beta^{-1}-1 +\sum_{i=0}^Z \mathcal{E}_{\infty}^{(i)}}\Bigr]},
\]
where $Z$ is distributed like the offspring of the Galton-Watson tree and is independent of all $\mathcal{E}_{\infty}^{(i)}$ which are independent copies with the same law as $\mathcal{E}_{\infty}$.
\end{theorem}

\begin{remark} In the case $\beta=1$, this formula allows one to recover a result from~\cite{LPP2} which states that $v(1)=E\Bigl[\frac{Z-1}{Z+1}\Bigr]$. Also, as previously mentioned this theorem implies Theorem~\ref{aid1}.\end{remark}

The central achievement of~\cite{Aidekon} is to describe the asymptotic distribution of the tree seen from the particle. In particular, A\"id\'ekon obtains a very explicit description of the environment seen from the particle up to orientation (i.e.~forgetting which path leads back to the root). This means that the average drift under the environment seen from the particle is explicit, and hence a formula for the speed can be found.

Even though the expression of the speed is not explicit enough to fully answer  the Question~\ref{conj_speed}, we are able to derive other interesting properties (Personal communication from Elie A\"id\'ekon) such as
\begin{itemize}
\item For a Galton-Watson tree with leaves, we have $\frac{v(\beta)-v(\beta_c)}{\beta-\beta_c}<-c$ with $c>0$, where $\beta_c$ was defined in Theorem~\ref{gw_lppspeed}.
\item For a Galton-Watson tree without leaves, the speed $v(\cdot)$ viewed as a function of $\beta$ has a continuous derivative on $(1,\infty)$.
 \end{itemize}

\subsubsection{The Einstein relation on Galton-Watson trees}\label{sect_arbre_einstein}

The Einstein relation is a general principle stemming from the fluctuation-dissipation theory, see~\cite{DD}. In general, it relates the diffusivity (the natural fluctuations of the system) to the response of the system (called the dissipation) to an external excitation, in this case, the addition of a small drift.  

\begin{remark}\label{bounds} In this section we will assume that the Galton-Watson trees do not have leaves and verify an exponential moment bound on the offspring distribution. These assumptions are probably unnecessary, a weak moment bound on the offspring distribution should be sufficient, but the results have not been rigorously  proved in another context. \end{remark}

In order to state the Einstein relation, we need to spend some time discussing the biased random walks in the null recurrent regime mentioned in Theorem~\ref{lycap}. A result from~\cite{PZ} shows that the $1/{\bf m}$-biased random walk on a Galton-Watson tree, which we call $X^{\text{rec}}_n$, verifies a quenched invariance principle (and thus an annealed one). More precisely, recalling the assumptions of Remark~\ref{bounds},
\begin{theorem}\label{PZ}
For ${\bf P}$-a.s.~every tree, the process $\Bigl\{\abs{X^{\text{rec}}_{\lfloor nt\rfloor}}/\sqrt{\sigma^2 n}\Bigr\}_{t\geq 0}$ converges in law, on the Skorohod space $D(\R^+,\R)$, to the absolute value of a standard Brownian motion. Moreover, the variance is explicit
\[
\sigma^2=\frac{m^2(m-1)}{\sum_{k=1}^{\infty} k^2 p_k-m}.
\]
\end{theorem}

In the context of biased random walks on Galton-Watson trees the Einstein relation was proved for continuous time versions of the walks. In that case, the invariance principle in Theorem~\ref{PZ} is preserved, although the limiting variance is multiplied by a factor $2$ where the factor $2$ is due to the speed up of the continuous-time walk relative to the discrete-time walk. This means that $X^{\text{rec}}_n$ rescaled by $\sqrt n$, behaves as a Brownian Motion of  diffusivity $\mathcal{D}=2\sigma^2$.

Introducing the notation $\overline{v}(\beta)$ for the speed of the continuous time $\beta$-biased random walk on the Galton-Watson tree, we are now able to state the Einstein relation proved in~\cite{BHOZ}. Recalling the assumptions of Remark~\ref{bounds}.
\begin{theorem}\label{BHOZ}
We have that
\[
\lim_{\alpha \to 0^+} \frac{\overline{v}(\exp(\alpha)/m)}{\alpha} =\frac{\mathcal{D}}2>0.
\]
\end{theorem}

In words let us explain what the previous relation means. Starting from the null-recurrent case and slightly increasing the bias will turn the system transient and, in fact, ballistic. The response to this excitation of the system (the speed) is directly linked to the diffusivity at equilibrium (in the null-recurrent case). In physics term, we say that the mobility (left-hand side) is related to the diffusivity (right-hand side). The constants in this relation are universal in the sense that they do not depend on the details of the system (here the precise law of the Galton-Watson tree).

 \begin{remark}
 Theorem~\ref{BHOZ} holds in a more general setting. It was actually proved on bi-infinite Galton-Watson trees, on which one can make sense of the speed even for bias lower that $1/{\bf m}$. In that case, the Einstein relation is proved to hold not only for the right-hand limit but also the left-hand one.
\end{remark}

\begin{remark} On trees, it is natural to expect that the Einstein relation holds in many related models. For example for biased random walks on Galton-Watson trees with leaves, on Galton-Watson trees with only moment bounds on the offspring distribution, or multi-type GaltonÐWatson trees as in~\cite{DSun}.
\end{remark}

As shown in the previous section, questions on the limiting velocity are extremely interesting. Nevertheless, they provide no information on the order of magnitude of $\abs{X_n}$ in the zero speed regime (except that it is sub-linear). We shall address this issue in the next section.

\subsection{Fluctuations for biased random walks on Galton-Watson trees.}

This section will mainly focus on fluctuations for biased random walks on Galton-Watson trees with leaves. Indeed, in the case without leaves the main questions are pretty much answered in~\cite{PZ}.
\begin{theorem}\label{PZ2}
Let us consider a Galton-Watson tree without leaves and exponential moment bounds on the offspring distribution and fix $\beta>1/m$. For ${\bf P}$-a.e. tree, the process $\Bigl\{\abs{X^{(\beta)}_{\lfloor nt\rfloor}-v(\beta)nt}/{\sqrt{\sigma^2 n}}\Bigr\}_{t\geq 0}$ converges in law to the absolute value of a standard Brownian motion.
\end{theorem}

With that being said, we will, henceforth, assume that our Galton-Watson tree has leaves.

Before stating the results on this model in Section~\ref{sect_arbre_result}, we shall first describe its link to trap models and discuss the lattice effect.

\subsubsection{Link to biased trap models}\label{sect_arbre_trap}

Let us explain how biased random walks on Galton-Watson trees can be viewed as a biased trap model on $\Z$. For this we need mainly to explain two facts
\begin{itemize}
\item How can we compare a walk on a tree and a trap model on $\Z$?
\item Where do polynomially  decaying tails appear?
\end{itemize}

The issues related to independence etc.~can be addressed in a similar manner to what was done for one-dimensional RWREs (see Section~\ref{sect_1dbtm}).

\vspace{0.4cm}

{\it How can we relate our model to a one-dimensional structure?}

\vspace{0.4cm}

Maybe the most surprising fact about biased trap models is their relative insensitivity to the underlying trap structure. Because of the exponential growth of the tree we know that there are, within distance $n^{\alpha}$ of the root, giant traps of height $H$ of polynomial order, and thus of waiting time that is at least exponential in $n^{\alpha}$. Why is it that we can ignore them?

The directed nature of the walk is the key factor here. Indeed, on a given realization, the directional transience will force the walk to move quickly through space and will only have the opportunity to visit a small fraction of the space. 

More specifically, for the $\beta$-biased random walk on a Galton-Watson tree, it can be shown that the number of different sites visited within time $n$ is of order $n$ (this is elementary if $\beta>1$, by comparing the walk on the backbone with the $\beta$-biased random walk on $\Z$). This vision should help to convince the reader of why the walk does, to some extent, behave in a similar manner on $\Z$ and on trees.

The way to make this rigorous is to use a tool called regeneration times. Directionally transient RWREs usually have a regenerative structure, which provides a key tool for their analysis, see~\cite{SZ}. Informally, we may define this structure as follows. Consider the first time at which the particle reaches a new maximum in the direction of the transience which is also
a minimum in this direction for the future trajectory of the random walk; we
call this time $\tau_1$. Even though $\tau_1$ is not a stopping time, it has the interesting
property of separating the past and the future of the random walk into two independent parts in an annealed sense.
Iterating this construction, we actually obtain a sequence $(\tau_k)_{k\in \N}$ of regeneration times, which separates the walk into independent blocks. This property can be extremely useful, because it reduces the problem of understanding a directionally transient random walk in random environment to a sum of i.i.d.~random variables.

In a word , using regeneration times, we can reduce our problem to a sum of i.i.d.~random variables.  The key issue is then to understand the tail of those random variables (i.e.~the time spent in regeneration) which essentially boils down to understanding the tail of the time spent in traps.

\begin{remark}\label{rem_diff_1d}
A quick word has to be said on the limitation of the comparison between biased trap models on $\Z$ and directionally transient RWREs. In general we believe that the annealed behavior should be similar, indeed the regeneration times bring the problem down to the understanding of a sum of i.i.d.~random variables. 

However the quenched behavior will differ in most cases when the underlying graph is $\Z$.

 Indeed, on $\Z$, the sequence of large traps the walker will visit is fixed before hand because there is only one path to infinity, this means that the time to reach level $n$ will be extremely dependent on the exact geometry of the big traps we have to meet. This explains why one dimensional RWREs have complicated quenched limiting behaviors, see Remark~\ref{extra_info}.
 
  In the case of trees, in order to reach level $n$, we have to choose one path among a large (exponential in $n$) number of possibilities. This means that the exact sequence of big traps the walk visits will be one chosen among a large number of possible sequences of big traps. In other words, the sequence of big traps visited is averaged over the environment, i.e.~chosen under the annealed measure. This is why we would expect, on trees and actually also on $\Z^d$ ($d\geq 2$), that the quenched behavior resembles the annealed one.
\end{remark}

{\it Polynomial tails}

\vspace{0.4cm}

The Harris-decomposition (see Section~\ref{sect_arbre_harris}) has given us a natural candidate for traps: it is a ${\bf h}$-Galton-Watson tree (which is finite). Given such a tree, we are interested in the time of an excursion of the biased random walk in this tree, which should correspond to the waiting times $\tau_i$ in the BTM (defined in Section~\ref{sect_btm}).

So we  generate an ${\bf h}$-Galton-Watson tree (we denote the corresponding law ${\bf P}_{\text{trap}}$), then we start a walk from the root and focus on the return time $T_{\text{trap}}$ to the root. 

Using the intuition from reversible Markov chain (more specifically the mean return formula), we will approximate $T_{\text{trap}} \approx \beta^H$, where $H$ is the height of the tree. More formally, $H =\max \{n \geq 0, Z_n >0 \}$ where $Z_n$ is the size of the $n$-th generation. By a result of~\cite{Heathcote}, we know that there exists $\alpha>0$ such that
\begin{equation}\label{eq_h}
{\bf P}_{\text{trap}}[H\geq n] \sim \alpha {\bf f}'(q)^n,
\end{equation}
using $T_{\text{trap}}\approx \beta^H$ (with abusive notations), we see that
\begin{equation}\label{fake_tail}
{\bf P}_{\text{trap}}[T_{\text{trap}} \geq n] \approx n^{- \ln \beta_c/\ln \beta},
\end{equation}
where $\beta_c=1/{\bf f}'(q)$. As expected in models with trapping, we see polynomial tails.

\begin{remark}Unsurprisingly, the expected trapping time becomes infinite as $\beta$ becomes larger than $\beta_c$ where $\beta_c$ appears in Theorem~\ref{gw_lppspeed}. This partially explains Theorem~\ref{gw_lppspeed}.
\end{remark}

After this small investigation, we are led to believe that our model should behave as a biased Bouchaud trap model with waiting times that behave like $\beta^{H}$ where $H$ is integer-valued and has geometric tail estimates (see~(\ref{eq_h})). There is, however, a subtle caveat in~(\ref{fake_tail}) that we swept under the rug and is critical for the limiting behavior. We shall discuss this issue in the coming section.

\subsubsection{The lattice effect}\label{sect_arbre_lattice}

In this section we are going to present the `lattice-effect', a phenomenon related to the limiting behavior of sums of i.i.d.~random variable. 

\vspace{0.4cm}

{\it A simple example}

\vspace{0.4cm}

Let us consider two sums of i.i.d.~random variables
\begin{itemize}
\item $S_1(n)=\sum_{i=1}^n \beta^{X^{(i)}}$ where the $X^{(i)}$ are i.i.d.~exponential random variables of parameter $\ln 2$,
\item $S_2(n)=\sum_{i=1}^n \beta^{Y^{(i)}}$ where the $Y^{(i)}$ are i.i.d.~geometric random variables of parameter $1/2$.
\end{itemize}

Although $X$ and $Y$ have very similar tails, the difference between them has a dramatic effect concerning the limiting behavior of $S_1(n)$ and $S_2(n)$. Indeed, a simple computation shows that
\[
{\bf P}[\beta^X\geq t] ={\bf P}[X\geq \ln t/\ln \beta] \sim t^{-\ln 2/\ln \beta},
\]
and, hence, $\beta^X$ belongs to the domain of attraction of a stable law (this expression is explained in the Appendix in Section~\ref{sect_def_domain}) by Theorem~\ref{sumiidgen}. But, on the other hand, keeping in mind that $Y$ is integer valued,
\begin{align*}
P[\beta^Y\geq t] &=P[Y\geq \ln t /\ln \beta] \\ = &P[Y\geq \lfloor \ln t /\ln \beta\rfloor] \sim 2^{- \lfloor \ln t /\ln \beta\rfloor} =t^{-\ln 2/\ln \beta}F(t),
\end{align*}
where $F(t)\in [1/2,1]$ but does not have slowly-varying tails (see the Appendix, definition~\ref{def_sec_varlente}, for a precise definition), which means that $\beta^Y$ does not belong to the domain of attraction of a stable law (this expression is explained in the Appendix in Section~\ref{sect_def_domain}) by Remark~\ref{no_dom_attract}. This will happen as soon as the distribution of $Y$ is concentrated on a lattice, hence the name lattice-effect.

This means that $S_1(n)$ and $S_2(n)$ do not have similar scaling limits and limiting behavior. 

\begin{remark}\label{no_scale} It is not possible to center and rescale $S_2(n)$ to obtain a scaling limit. \end{remark}

\begin{remark} The lattice-effect arises from the lattice nature of an element in the environment. It is not linked to the choice of a discrete time Markov chain. \end{remark}

\begin{remark} The hypothesis that the distribution of $\log \rho_0$ be non-lattice (that is, the $\Z$-linear span of $\log \text{supp}(\rho_0)$ is dense in $\R$) in Theorem~\ref{KKS} is exactly what is needed to allow for scaling limits to exist. \end{remark}

\vspace{0.4cm}

{\it A heuristic for this phenomenon}

\vspace{0.4cm}

Let us explain why Remark~\ref{no_scale} is not such a surprising fact after all. 

Consider the case of $S_2(n)=\sum_{i=1}^n \beta^{Y^{(i)}}$ where the $Y^{(i)}$ are i.i.d.~geometric random variables of parameter $1/2$. It typically takes $i_0(n)$ indices to observe the first $Y^{(i)}$ equal to $n$. We can notice that
\begin{enumerate}
\item the first $Y^{(i)}$ that equals $n+1$ typically occurs around the index $2i_0(n)$,
\item right before this we have witnessed only $Y^{(i)}$s  lower than $n$, and usually only 2 of them have value $n$.
\end{enumerate}

Taking $\beta$ to be very large, $S_2(n)$ should be dominated by the terms with the largest values of $Y^{(i)}$ so it seems fair to say that, for a small $\epsilon>0$, we have
\[
S_2(2(1-\epsilon) i_0(n))\approx 2 \times \beta^n \text{ and } S_2(2(1+\epsilon) i_0(n))\approx \beta \times \beta^n.
\]

It should be clear from the previous equation that any sort of scaling will be difficult, since in a very short period of time the sum increases enormously. Those large changes will happen with an exponential periodicity.

\begin{remark} It is possible to obtain scaling limits for $S_2(n)$ along exponentially growing subsequences (of the form $\lambda 2^n$ in our example) towards infinitely divisible distributions. This is a consequence of Theorem IV.6 (p.77) in~\cite{Petrov}.
\end{remark}

\subsubsection{Fluctuations for the biased random walk on the Galton-Watson tree with leaves in the sub-ballistic regime}\label{sect_arbre_result}

In Section~\ref{sect_arbre_trap}, we argued that biased random walks on Galton-Watson trees with leaves should behave like $\sum_{i=0}^n \beta^{H^{(i)}}$ where $H^{(i)}$ are i.i.d.~integer valued with a tail given by $P[H\geq n] \sim \alpha {\bf f}'(q)^n$, as we recall from~(\ref{eq_h}).

In Section~\ref{sect_arbre_lattice}, we then discussed how this setting differs from that of the traditional biased BTM because the distribution of $H$ is concentrated on a lattice. Because of this, we should expect scaling limits not to exist, except maybe along exponentially growing subsequences. 


\vspace{0.4cm}

{\it Results on the biased random walk on the Galton-Watson tree with leaves in the sub-ballistic regime}

\vspace{0.4cm}

We recall that the model of biased random walks on Galton-Watson trees was defined in Section~\ref{def_model} and that all necessary notations can be found in that part of the notes.

Throughout this section we assume a moment bound on the offspring distribution namely $E[Z^2] <\infty$. Before moving onto more difficult questions, let us notice that it was shown in~\cite{FG} that an annealed central limit theorem holds for the $\beta$-biased random walk on the Galton-Watson tree with leaves provided $\beta<\beta_c^{1/2}$ where $\beta_c$ appeared in Theorem~\ref{gw_lppspeed}. This corresponds to the cases $\alpha>2$ in the coming equation~(\ref{alpha_tree}).

Assume that $\beta > \beta_c=1/ {\bf f}'(q) $ and introduce 
\begin{equation}\label{alpha_tree}
\alpha:= \frac{-\ln {\bf f}'(q)}{\ln \beta}=\frac{\ln \beta_c}{\ln \beta} < 1,
\end{equation}
which is reminiscent of~(\ref{fake_tail}). We denote $\Delta_n$ be the hitting time of the $n$-th level:
\[
\Delta_n=\inf\{i\geq 0: \abs{X_i}=n\}.
\]

First, we see that $n^{1/\alpha}$ is the correct scaling. Indeed, we have the following two results from~\cite{FG},
\begin{theorem} For any $\beta>\beta_c$.
\begin{enumerate}
\item The laws of $(\Delta_n/n^{1/\alpha})_{n\geq 0}$ under $\PR$ are tight.
\item The laws of $(\abs{X_n}/n^{\alpha})_{n\geq 0}$ under $\PR$ are tight.
\item We have $\lim_{n \to \infty} \frac{\ln \abs{X_n}}{\ln n}=\alpha,\ \PR-\text{a.s.}$
\end{enumerate}
\end{theorem}

But, as expected from our heuristic, the following was proved in~\cite{FG}
\begin{theorem}
For $\beta$ large enough, the sequence $(\Delta_n/n^{1/\alpha})_{n\geq 0}$ does not converge in distribution. 
\end{theorem}

\begin{remark} We believe that the previous theorem should be true for any $\beta>\beta_c$.\end{remark}

Furthermore,
\begin{theorem}
\label{subsequ}
For any $\lambda>0$, denoting $n_{\lambda}(k)=\lfloor \lambda {\bf f}'(q)^{-k} \rfloor$, we have
\[
\frac{ \Delta_{n_{\lambda}(k)}} {n_{\lambda}(k)^{1/\alpha}} \xrightarrow{(d)} Y_\lambda
\]
where the random variable $Y_\lambda$ has an infinitely divisible law $\mu_\lambda$.
\end{theorem}

We now describe the limit laws $\mu_\lambda$ which has a particular form. There are constants $C_1,C_2$ such that
we have 
\[
Y_\lambda = (C_1\lambda)^{1/\alpha}\,\,  \widetilde{Y}_{(C_1C_2\lambda)^{1/\alpha}},
\]
where 
\[
\widetilde{Y}_\lambda \text { has the law } \IG(d_\lambda,0,\LR_\lambda).
\]

The infinitely divisible law 
$\IG(d_{\lambda},0,\LR_\lambda)$
is given by its L\'evy representation (see \cite{Petrov}, p. 32).
More precisely, the characteristic function of $\IG(d_{\lambda},0,\LR_\lambda)$ can be written in the form
\[
\ES\left [e^{it\widetilde{Y}_\lambda}\right] = \int e^{itx}\IG(d_{\lambda},0,\LR_\lambda) (dx) = 
\exp\left(id_\lambda t + \int\limits_0^\infty \left(e^{itx}-1 - \frac{itx}{1+x^2}\right)d\LR_\lambda(x)\right)
\]
where $d_\lambda$ is a real constant and $\LR_\lambda$ a real function which is non-decreasing on the interval $(0, \infty)$ and satisfies $\LR_\lambda(x) \to 0$ for $x \to \infty$ and $\int \limits_0^a x^2 d\LR_\lambda(x) < \infty$ for every $a > 0$. Comparing to the general representation formula in \cite{Petrov}, p. 32, we here have that the gaussian part vanishes and $\LR_\lambda(x) =0$ for $x < 0$. The function $\LR_\lambda$ is called the L\'evy spectral function.
Note that $\LR_\lambda$ is not a L\'evy-Khintchine spectral function.

Hence, as can be seen biased random walks on Galton-Watson trees with leaves do not belong to the same universality class as the biased BTM or one-dimensional random walks. Nevertheless, with a very small change in the definition of the model things would be different. This is illustrated by randomly biased random walks on Galton-Watson tree with leaves.

\vspace{0.4cm}

{\it Randomly biased random walks on Galton-Watson tree with leaves} 

\vspace{0.4cm}

Let us now present a result contained in two articles~\cite{H} and~\cite{BH}. Instead of taking a constant bias $\beta$,  we could choose to assign to each edge in the tree a random bias chosen according to some law $\nu$ with support in $(1+\epsilon, M)$ for some $\epsilon >0$ and $M<\infty$. For a precise definition, see~\cite{H} or consult Section~\ref{sect_arbre_random} for related models. 

Our entire heuristic on the bias provoking slowdown in leaves should remain true, but, the randomization of the bias could be a way to get rid of lattice-effects. And, indeed, see~\cite{H}
\begin{theorem}\label{alan}
Under some mild integrability condition on the offspring distribution (apart from $E[Z]>1$ and $p_0>0$).

Assuming also that the support of $\nu \circ \log^{-1}$ is non-lattice (that is, the $\Z$-linear span of $\log \text{supp}(\nu)$ is dense in $\R$), we have the following:  define $\alpha$ to be such that
\[
\int_0^{\infty} y^{\alpha} \nu(dy)=\frac 1 {{\bf f}'(q)},
\]
where $q$ is the extinction probability of the Galton-Watson tree and ${\bf f}$ its generating function. Then, if $\alpha<1$, we have the annealed convergences
\[
\frac{\Delta_n}{n^{1/\alpha}} \xrightarrow{(d)} C  \mathcal{S}_{\alpha}^{ca},
\]
where $\mathcal{S}_{\alpha}^{ca},$ is a completely asymmetric $\alpha$-stable law (see the Appendix, Section~\ref{def_stable_law_section}, for a more precise definition) and
\[
\frac{\abs{X_n}}{n^{\alpha}} \xrightarrow{(d)} C^{-\alpha} ( \mathcal{S}_{\alpha}^{ca})^{-\alpha}.
\]

Moreover, if $\alpha>1$ then the motion is ballistic.
\end{theorem}

\subsubsection{A more detailed analysis of the time spent in traps}\label{sect_approx}

We will give here a sketch on how to analyze the time spent in traps in the case of biased random walks on Galton-Watson trees with leaves. The method presented here can, up to minor modifications, be adapted to different contexts. We will, from time to time, give arguments in our sketch that are slightly too general for the model that we are considering, this is to give a flavor of potential generalizations.

The structure of large traps (subcritical Galton-Watson trees) is well understood, see~\cite{Geiger}. It can be seen as a long one dimensional line with comparatively small sub-traps (sub-trees) hanging off it.  On Figure~11 we have drawn a typical big trap emphasizing the key elements: the top, the height and the bottom.

\begin{figure}
\centering 
\epsfig{file=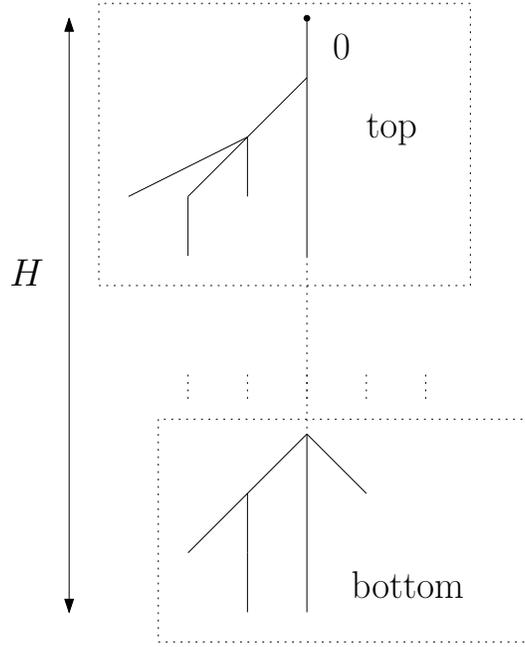, scale=0.7}
\caption{The critical parts for trapping of a large subcritical Galton-Watson tree}
\end{figure}

Let us now take a look at the law of the time of an excursion in a typical trap.
\begin{enumerate}
\item Either the walk does not hit the bottom before coming back to the top, this has a probability $p_{\text{bottom}}$ that is well approximated by information at the top of the tree. In this case, the walk only has access to small sub-traps and the return time to the origin is small.
\item Otherwise the walk hits the bottom of the trap and the time to descend is quick. From there the walk will start doing excursions from the bottom. The number of such excursions is a geometric random variable with parameter $p_{\text{top}}$ where we can write $p_{\text{top}}=\tilde{p}_{\text{top}}\beta^{-H}$ where $\tilde{p}_{\text{top}}$ is of order one and essentially measurable with respect to information from the top of the tree ($p_{\text{top}}$ and $p_{\text{bottom}}$ can be linked by reversibility). All those excursions are i.i.d.~and their duration $T^{(i)}_{\text{exc}}$ has a mean that can be well estimated (using the mean return time formula) with information from the bottom of the tree. Eventually, the walk will go back up and, conditionally on not returning to the bottom, that trip up goes fast.
\end{enumerate}

This long description leads to the following approximation (using the notations $\text{Bern}(p)$ for a Bernoulli random variable of parameter $p$ and $\text{Geom}(p)$ for a Geometric random variable of parameter $p$)
\begin{align*}
T_{\text{trap}} & =\1{\text{Bern}(p_{\text{bottom}})=1} \sum_{i=1}^{\text{Geom}(\beta^{-H} \tilde{p}_{\text{top}})}  T^{(i)}_{\text{exc}} \\
                         &\approx \1{\text{Bern}(p_{\text{bottom}})=1} \text{Geom}(\beta^{-H} \tilde{p}_{\text{top}})  E^{\omega}[T^{(i)}_{\text{exc}}] \\
                          &\approx  \1{\text{Bern}(p_{\text{bottom}})=1} \beta^{H} \tilde{p}_{\text{top}}^{-1} E^{\omega}[T^{(i)}_{\text{exc}}] {\bf e},
                          \end{align*}
where ${\bf e}$ is an exponential random variable of mean $1$. Here, we used the law of large numbers in the second line and the approximation of geometric random variable of small parameters by exponentials. For completeness let us specify the measurability issues
\begin{itemize}
\item $\text{Bern}(p_{\text{bottom}})$ depends on the randomness of the walk and its parameter $p_{\text{bottom}}$ on the top of tree,
\item $\tilde{p}_{\text{top}}^{-1}$ depends on the top of the tree,
\item $E^{\omega}[T^{(i)}_{\text{exc}}]$ depends on the bottom of the tree
\item ${\bf e}$ depends on the randomness of the walk.
\end{itemize}

The top, the bottom and the height  of the trap  are asymptotically independent (and obviously independent of the walk). Hence, we can deduce the tail of $T_{\text{trap}}$ from the tail of $H$ (after proving some pretty weak moment bounds on the other variables involved).

The last step needed to find the total time spent in a given trap is to understand the number of entries in a deep trap. This argument typically involves understanding something about the environment seen from the particle at the top of the trap. This matter can be rather subtle and we do not pursue it here.

\begin{remark}\label{1d_trapex} In the model of one-dimensional RWREs, the bottom of the trap is naturally the bottom of the valley (see Section~\ref{sect_1dbtm} for the associated vocabulary) and the top of the trap corresponds to the right-hand side of the valley  after which the potential decreases.
\end{remark} 

\subsection{Other models of random walks in random environments on trees}\label{sect_arbre_random}

In this section, we will present randomly biased random walks on Galton-Watson trees without leaves.

\subsubsection{The model}\label{model_aidekon}

We construct a supercritical Galton-Watson tree without leaves, i.e.~$p_0=0$, (similar models can be considered with leaves but we choose to ignore this) whose mean we denote ${\bf m}$. Furthermore consider a sequence of i.i.d.~random variables $A^{\infty}:=(A^{(i)})_{i\geq 1}$ taking values in $\R^+$ with a common law denoted $A$. We will mark the vertices of the tree using those random variables.

At the root, we pick  random variables $(A_0^{(i)}$, with $i \leq Z(0))$ (with the same law as $A^{\infty}$), assigning a coefficient $A_0^{(i)}$ to $v_i$  where $v_i$ is one of the $Z(0)$ descendants of the root.  For each vertex $u$ of the $n$-th generation, we pick independently a random vector $A_u^{(i)}$ (with the same law as $A^{\infty}$), with $i \leq Z(u)$, and assign a coefficient $A_u^{(i)}$ to the edge $v_i$  where $v_i$ is one of the $Z(u)$ descendants of $u$. 

In this marked environment, for a vertex $u$ that has $k$ children $v_1,\ldots, v_k$ and parent $\overleftarrow{u}$, then the transition probabilities for our Markov chain are given by
\begin{enumerate}
\item $P^{\omega}[X_{n+1}=\overleftarrow{u}|X_n=u]=\frac 1 {1+\sum_{j=1}^k A_u^{(j)}}$, 
\item $P^{\omega}[X_{n+1}= v_i |X_n=u] = \frac {A_u^{(i)}} {1+\sum_{j=1}^k A_u^{(j)}}$, for $1\leq i\leq k$.
\end{enumerate}

Even though the rules are not well defined at the origin this can be dealt with by adding artificially a vertex $\overleftarrow{0}$, with the transition rule $P^{\omega}[X_{n+1}=0|X_n=\overleftarrow{0}]=1$.

We will assume in the sequel that $\text{supp}(A)\in (\epsilon,M)$ with $\epsilon>0$ and $M<\infty$. 

\subsubsection{Results}

The most basic result on this model goes back to~\cite{LyonsPemantle}.
\begin{theorem}[Lyons, Pemantle - 1992]
We have the following criterion
\begin{enumerate}
\item if $\inf_{t\in [0,1]} {\bf E}[A^t]>1/m$ then $X_n$ is transient,
\item if $\inf_{t\in [0,1]} {\bf E}[A^t]<1/m$ then $X_n$ is recurrent.
\end{enumerate}
\end{theorem}

Let us now present a result of A\"id\'ekon~\cite{Aidekon1}.  For this, we introduce $\alpha=\alpha_1+\alpha_2$ (the usefulness of this notation will become clear in the next section) with
\begin{equation}\label{aid_alpha}
\alpha_1:=\text{Leb}\Bigl\{t\geq 0,\ {\bf E}[A^t] \leq 1/p_1\Bigr\} \text{ and } \alpha_2:=\text{Leb}\Bigl\{t\leq 0,\ {\bf E}[A^t] \leq 1/p_1\Bigr\},
\end{equation}
where we set $\alpha:=\alpha_1:=\alpha_2:=\infty$ if $p_1=0$ (we recall that $p_1$ is the probability to have one offspring). Then
\begin{theorem}
In the transient case, we have
\begin{itemize}
\item if $\alpha>1$, then $\lim_{n\to \infty} \abs{X_n}/n=v>0$,
\item if $\alpha<1$, then $\lim_{n\to \infty} \abs{X_n}/n=0$, and
\[
\lim_{n\to \infty} \frac{\ln (\abs{X_n})}{\ln n} =\alpha, \qquad \PR\text{-a.s.}.
\]
\end{itemize}
\end{theorem}

The result certainly indicates that this model should exhibit trapping. This is natural, indeed since the bias is allowed to be upwards, we do not require leaves to create traps.

\subsubsection{Relation to trapping}

The trapping in this model is different from the examples we previously encountered, let us start by describing the mechanisms slowing the walk down before presenting related open questions.

\vspace{0.4cm}

{\it Link to trapping}

\vspace{0.4cm}

 In the work~\cite{Aidekon1}, it is shown that a slowdown (provoking sub-ballistic behavior) can only occur in long pipes  (sections of the tree with only one child) inside the tree. This explains the importance of $p_1$ in~(\ref{aid_alpha}). 
 
 With this information in mind, we can explain what occurs (we stress that the following is conjectured and not proved). In those pipes, there is a potential defined in a manner very similar to the way it was done in one dimension (see~(\ref{potential}) replace the $\rho$s by $A$s). For a trap to form inside a pipe, we need for the potential in that pipe to form a valley. The cost of forming a pipe with such a valley of height $H\geq n$ comes from two factors
 \begin{itemize}
 \item we need at first a descending slope of height $n$ (whilst forcing the tree to have only one offspring), this should cost roughly $\exp(-\alpha_1 n)$,
\item we then need to form an ascending slope of height $n$ (whilst forcing the tree to have only one offspring), this should cost roughly $\exp(-\alpha_2 n)$.
\end{itemize}

This means that the total cost of the trap of height $H\geq n$ is roughly $\exp(-\alpha n)$. Considering the resemblance with~(\ref{1d_height}), we think this sketch should be sufficient to convince the reader that the link to trapping models is pretty explicit.

\subsection{Open problems}

Let us present some of the questions that have not yet been addressed in the literature and that could be accessible using the approach of trap models.
\subsubsection{Summary of results and remaining open questions}

Many of the following questions first appeared in~\cite{LPP1}, where many more issues concerning random walks on trees are considered. 

\vspace{0.4cm}

{\it Galton-Watson trees without leaves}

\vspace{0.4cm}

Let us first show Figure~12 which depicts of the conjectured behavior of the speed. 

\begin{figure}
\centering 
\epsfig{file=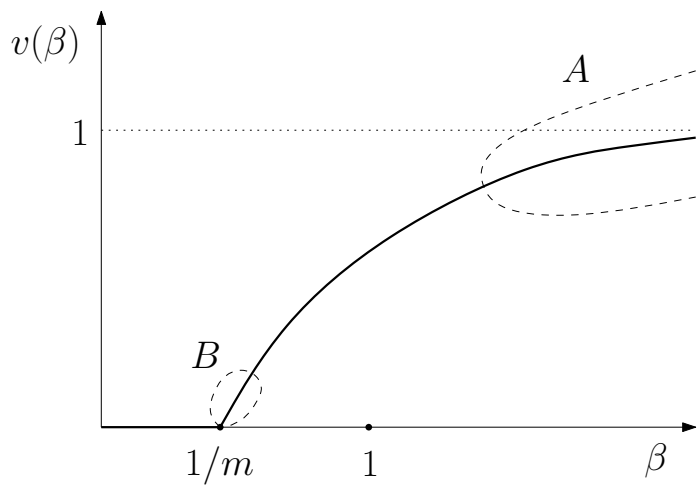, scale=0.7}
\caption{The conjecture for the behavior in terms of $\beta$ of the speed of a biased random walk on a Galton-Watson tree without leaves}
\end{figure}

The main questions on the speed as a function of $\beta$ are the following.
\begin{question}
 Is the speed monotonic in the bias? 
 \end{question}
 
This is answered in the zone A (by Theorem~\ref{BAFS} or Theorem~\ref{aid1}) and it would be in zone B (by Theorem~\ref{BHOZ}) provided we knew that $v(\cdot)$ as a function of $\beta$ has a continuous derivative in a neighborhood of $1/{\bf m}$.

\begin{question}
 What are the regularity properties of the speed? 
  \end{question}
  
 This is mainly open. We only know that for a Galton-Watson tree without leaves, the speed $v(\cdot)$ viewed as a function of $\beta$ has a continuous derivative on $(1,\infty)$ (this is a consequence of  Theorem~\ref{aidekon}).

\begin{question}
  What is the right derivative of $v(\cdot)$ at $1/{\bf m}$?
\end{question}

 This is answered by Theorem~\ref{BHOZ}.

Another natural question concerning the speed as a function of the environment is
\begin{question}
 Is the speed monotonic in the environment? 
\end{question}

This is partially answered in Theorem~\ref{NYU} and Remark~\ref{FL}.

\vspace{0.4cm}

{\it Galton-Watson trees with leaves}

\vspace{0.4cm}

The expected behavior of the speed in this context is presented in Figure~13.

\begin{figure}
\centering 
\epsfig{file=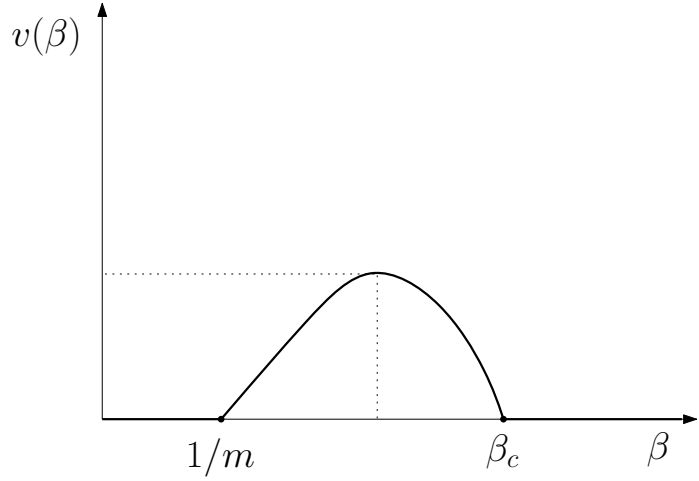, scale=0.7}
\caption{The conjecture for the behavior in terms of $\beta$ of the speed of a biased random walk on a Galton-Watson tree with leaves}
\end{figure}

The main questions on the speed as a function of $\beta$ are the following
\begin{question}
 Is the speed unimodal? 
 \end{question}
 
 To the best of our knowledge, this is completely open.
 
 \begin{question}
  What can we say about the right (resp.~left) derivative of the speed at $1/{\bf m}$ (resp.~$\beta_c$)?
  \end{question}
  
   The value of the first derivative should follow from simple modifications of the work in Theorem~\ref{BHOZ}, the derivative at $\beta_c$ provided it exists would be negative  by a consequence of Theorem~\ref{aidekon}.
   
   \begin{question}
 If the speed is unimodal, where is the maximum of the speed located, and what is this maximal value? 
 \end{question}
 
 As far as we know, nothing is known about this.

\begin{question}
 What are the regularity properties of the speed? 
  \end{question}
  
Once again nothing is known about this.

\subsubsection{Scaling limits}

Essentially, we wish to prove something similar to the picture provided by the biased BTM (stable limit laws, stable subordinators and aging). In the lattice-case of course, the results should be modified to adjust for the fact that our trapping times are not in the domain of attraction of stable laws.

\vspace{0.4cm}

{\it The central limit theorem}

\begin{conjecture}
A quenched central limit theorem should hold for the $\beta$-biased random walk on the Galton-Watson tree with leaves provided $\beta<\sqrt {\beta_c}$. 
\end{conjecture}

 This would extend the result mentioned at the beginning of Section~\ref{sect_arbre_result}. It should be provable with a reasoning resembling that of Remark~\ref{rem_diff_1d}.

\begin{conjecture}
In the case with randomly biased random walk, a quenched central limit theorem should also hold for $\alpha>2$ (where $\alpha$ was defined in Theorem~\ref{alan}).
\end{conjecture}

{\it Stable fluctuations}

\begin{conjecture}
In the case of a $\beta$-biased random walk with an exponent $\alpha \in (1,2)$ (defined at~(\ref{alpha_tree})), we would expect that there does not exist a way to recenter and rescale the walk to obtain scaling limits. Nevertheless, by re-centering $\Delta_n$ by the inverse of the speed and rescaling by $n^{1/\alpha}$, we should obtain a tight sequence that converges along exponential subsequences to an infinitely divisible distribution. 
\end{conjecture}

\begin{conjecture}
In the case of randomly biased random walk, the picture of the biased BTM should be correct. For $\alpha\in (1,2)$, (as defined in Theorem~\ref{alan}) a re-centerring of $\abs{X_n}$ by the speed and a rescaling by $n^{1/\alpha}$ should lead to stable scaling limits. At the process level, we would expect to see an $\alpha$-stable asymmetric L\'evy process in the limit.
\end{conjecture}

{\it Aging}

\begin{conjecture}
In the case of randomly biased random walks for $\alpha\in (0,1)$, an aging similar to the result in Theorem~\ref{aging_esz} should be true.
\end{conjecture}

{\it Cases $\alpha=1$ and $\alpha=2$}

\begin{question}
Can one prove more results in all the models previously mentioned when $\alpha=1$ or $\alpha=2$ in random walks in random environments?
\end{question}

To the best of our knowledge only for the RWRE on $\Z$ do we have results for $\alpha=1$ and $\alpha=2$ (see~\cite{KKS}). Using the approach from trap models only the case $\alpha=1$ has been studied for one-dimensional RWREs in~\cite{ESTZ}.

\subsubsection{Randomly biased random walks on Galton-Watson trees without leaves}

As far as the model of randomly biased random walks on Galton-Watson trees without leaves (presented in Section~\ref{model_aidekon}) is concerned, we can make the following conjectures.

\begin{conjecture} Adding a non-lattice condition on $\log A$ should allow one to prove that this model  has stable limit laws, stable  subordinators as limiting processes and verifies the aging property when $\alpha \in (0,1)$. 
\end{conjecture}

\begin{conjecture} It would be interesting to investigate the case $\alpha \in (1,2)$ and see if stable fluctuation occurs, assuming non-latticity. This should hold only if a slowdown strong enough to break the CLT cannot occur outside of pipes. 
\end{conjecture}

\section{Biased random walks on supercritical percolation clusters}\label{sec_percomain}

As mentioned in the introduction (see Section~\ref{sect_plan}) the model of biased random walks on supercritical percolation clusters is one of the motivations that sprung a lot of the work that has been done on trees. In this section, we present the current state of the research on this model mentioning that research is ongoing and that we expect new results to come out within a few years.

\subsection{Model}\label{perco}

Firstly, we describe the environment. We denote by $E(\Z^d)$ the edges of the nearest-neighbor lattice $\Z^d$ for some $d\geq 2$. We fix $p \in (0,1)$ and perform a Bernoulli bond-percolation by picking a random configuration $\omega \in \Omega:=\{0,1\}^{E(\Z^d)}$ where each edge $e$ has probability $p$ of verifying $\omega(e)=1$, independently of the assignations made to all the other edges. We introduce the corresponding measure 
\[
P_{p}= (p \delta_1 + (1-p) \delta_0)^{\otimes E(\Z^d)}.
\]

An edge $e$ will be called open in the configuration $\omega$ if $\omega(e)=1$. The remaining edges will be called closed.  This naturally induces a subgraph of $\Z^d$ which will also be denoted by $\omega$, and it  yields a partition of $\Z^d$ into connected components, called open clusters, and isolated vertices. 

It is classical in percolation that for $p>p_c(d)$, where $p_c(d)\in (0,1)$ denotes the critical percolation probability of $\Z^d$ (see~\cite{Grimmett} Theorem~1.10 and Theorem~8.1), there exists a unique infinite open cluster $K_{\infty}(\omega)$, $P_p$ almost surely. Moreover, the following event has positive $P_p$-probability: 
\[
\mathcal{I}= \Big\{ \text{there is a unique infinite cluster $K_{\infty}(\omega)$ and it contains $0$} \Big\}.
\]

 We further define
\[
{\bf P}_p[~\cdot~] = P_{p}[~\cdot\mid \mathcal{I}].
\]

The bias $\ell=\lambda \vec \ell$ depends on two parameters: the strength $\lambda >0$, and the bias direction $\vec \ell \in S^{d-1}$ which lies in the unit sphere with respect to the Euclidean metric of $\R^d$. Given a configuration $\omega \in \Omega$, we consider the reversible Markov chain $X_n$ on $\Z^d$ with law $P^{\omega}$,  whose transition probabilities $p^{\omega}(x,y)$ for $x,y\in \Z^d$ are defined by
\begin{enumerate}
\item $X_0=0$, $P^{\omega}$-a.s.,
\item $p^{\omega}(x,x)=1$, if $x$ has no neighbor in $\omega$, and
\item $\displaystyle{p^{\omega}(x,y) =\frac{c^{\omega}(x,y)}{\sum_{z \sim x} c^{\omega}(x,z)}}$,
\end{enumerate}
where $x\sim y$ means that $x$ and $y$ are adjacent in $\Z^d$.  Here we set
\[
\text{for all $x,y \in \Z^d$,} \qquad \ c^{\omega}(x,y)=\begin{cases}
                                            e^{(y+x)\cdot\ell} & \text{ if } x\sim y \text{ and } \omega(\{x,y\})=1, \\
                                            0           & \text{ otherwise.}\end{cases} 
\]                                           

The random variable $c^{\omega}(x,y)$ is called the conductance of the edge $e=[x,y]$ in the configuration $\omega$, a notation which is natural in light of the relation between reversible Markov chains and electrical networks,
for a presentation of which the reader may consult
 ~\cite{DoyleSnell} or ~\cite{LP}. The above Markov chain is reversible with respect to the  invariant measure given by 
\[
\pi^{\omega}(x)=\sum_{y \sim x} c^{\omega}(x,y).
\]  

Finally, we define the annealed law of the biased random walk on the infinite percolation cluster by the semi-direct product $\PR = {\bf P}_p  \times  P^{\omega}[\,\cdot\,]$.

\begin{remark}\label{BGP}
A similar model was considered in~\cite{BGP}, where the conductances are $c^{\omega}(x,y)=\beta^{\max (x\cdot e_1,y\cdot e_1)} \1{\omega(\{x,y\})=1}$ with $\beta>1$. This model is closer to our counterparts considered on trees but do not allow for the bias to be in any directions.
\end{remark}

\subsection{Results}

Physicists had long conjectured the existence of sharp phase transition from positive to zero speed, see~\cite{BD},~\cite{Dhar} and~\cite{DS}. The mathematical proof of this fact arrived in two steps

\subsubsection{Phase transition of speed}

In order to state the results, let us recall some results proved in~\cite{Sznitman}. In all cases, i.e., for $d\geq2$ and $p > p_c$, the walk is directionally transient in the direction $\vec{\ell}$ (see Theorem 1.2 in~\cite{Sznitman}):
 \begin{equation}
 \label{def_dir_trans}
 \lim_{n\to \infty} X_n \cdot \vec{\ell} =\infty, \qquad \PR-\text{almost surely,}
 \end{equation}
 and verifies  the law of large numbers  (see Theorem 3.4 in~\cite{Sznitman}),
\[
\lim_{n\to \infty} \frac{X_n} n =v , \qquad \PR-\text{almost surely},
\]
where $v\in \R^d$ is a constant vector. However, the main result in~\cite{Sznitman} was to prove that there exists a phase transition from positive to zero speed.
\begin{theorem}
There exist $0<\lambda_1\leq \lambda_2<\infty$ such that
\begin{itemize}
\item if $\lambda<\lambda_1$, then $v\cdot \vec{\ell}>0$,
\item if $\lambda >\lambda_2$, then $v=\vec{0}$.
\end{itemize}
\end{theorem}

\begin{remark} In~\cite{BGP}, the authors show a very similar statement in their model  (mentioned in Remark~\ref{BGP}) in $\Z^2$. \end{remark}

A big remaining conjecture was to prove that the previous phase transition was sharp, i.e.~that we can choose $\lambda_1=\lambda_2$ in the previous theorem. This seems natural in light of our results on trees and trap models. Actually, before any explicit link between traps models and RWREs had been proved, A-S. Sznitman mentioned in the survey \cite{AlainHouches} that the slowdown effects in this model  were similar to those responsible for aging in BTM.

\subsubsection{Sharpness of the phase transition}
 
We now introduce the backtrack function of $x \in \Z^d$, which will be fundamental for gauging the extent of the slowdown effect on the walk: 
\begin{equation}
\label{def_backtrack}
\mathcal{B}\mathcal{K}(x)=\begin{cases} 0 & \text{ if } x\notin  K_{\infty}\\
                                              \displaystyle{\min_{(p_x(i))_{i\geq 0}\in \mathcal{P}_x}}\max_{i\geq 0} (x-p_x(i))\cdot \vec{\ell}  & \text{otherwise,}
\end{cases}
\end{equation}
where $\mathcal{P}_x$ is the set of all infinite open vertex-self-avoiding paths starting at~$x$. As Figure~14 indicates, connected regions where $\backtrack$ is positive may be considered to be traps: indeed, from points in such regions, the walk must move in the direction opposed to the bias in order to escape the region. The height of a trap may be considered to be the maximal value of $\backtrack (x)$ attained by the vertices $x$ in the trap. It was proved in~\cite{FH} that
\begin{proposition}
\label{backtrack_exponent}
For $d\geq 2$, $p > p_c$ and $\vec{\ell} \in S^{d-1}$, there exists $\hatxi(p,\vec{\ell},d)\in (0,\infty)$ such that 
\[
\lim_{n \to \infty}  n^{-1} \log {\bf P}_p[ \mathcal{B}\mathcal{K}(0)> n  ] = -\hatxi(p,\vec{\ell},d).
\]
\end{proposition}

\begin{figure}\label{defBK}
\centering\epsfig{file=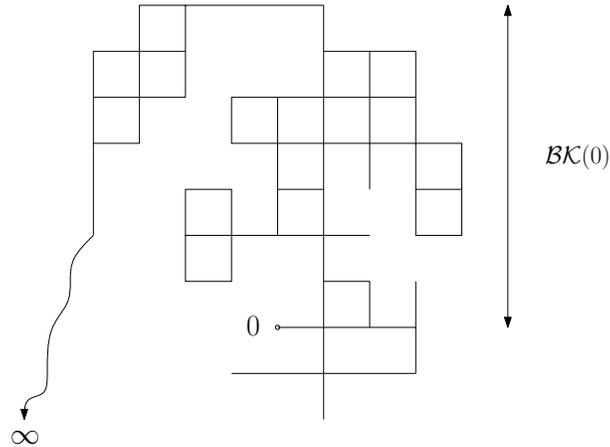, width=8cm}
\caption{Illustrating the value of $\mathcal{B}\mathcal{K}(0)$ in the case of an axial bias $\vec \ell$ pointing downwards.}
\end{figure}

This allows us to define the exponent $\alpha \in (0,\infty)$:
\begin{equation}
\label{def_gamma}
\alpha= \frac {\hatxi}{2\lambda}.
\end{equation}

The main result of~\cite{FH} gives a sharp transition from  a ballistic to a sub-ballistic regime. Moreover, the critical value~$\lambda_c$ of the bias is explicitly computed.
\begin{theorem}
\label{the_theorem}
For $d\geq 2$ and $p > p_c$, set $\lambda_c=\hatxi(p,\vec{\ell},d)/2$. We have that
\begin{enumerate}
\item if $\lambda <\lambda_c$, or, equivalently, $\alpha>1$, then $v\cdot \vec{\ell} >0$,
\item  if $\lambda >\lambda_c$, (or $\alpha<1$), then $v=\vec{0}$.
\end{enumerate}
\end{theorem}

Another by-product of the techniques of~\cite{FH} is that
\begin{theorem}
\label{the_theorem2}
Let $d \geq 2$ and $p > p_c$.  
Assume that $\alpha >2$.

 The $D(\R_+,\R^d)$-valued processes $B_{\cdot}^n = \frac 1 {\sqrt n}(X_{[\cdot n]}-[\cdot n] v)$ converge under $\PR$ to a Brownian motion with non-degenerate covariance matrix, where $D(\R_+,\R^d)$ denotes the space of right continuous $\R^d$-valued functions on $\R_+$ with left limits, endowed with the Skorohod topology, c.f.~Chapter 3 of \cite{EK}.
\end{theorem} 

Finally, in the sub-ballistic regime, the  polynomial order of the magnitude of the walk's displacement is given in~\cite{FH}.
\begin{theorem}
\label{the_theorem3}
Set $\Delta_n = \inf \big\{ m \in \N: X_m \cdot \vec{\ell} \geq n \big\}$. 

Let $d\geq 2$ and $p > p_c$.  
If $\alpha \leq 1$ then
\[
\lim \frac{\log \Delta_n}{\log n} = 1/\alpha, \qquad \PR-\text{almost surely},
\]
and
\[
\lim \frac{\log X_n \cdot \vec{\ell}}{\log n} =\alpha, \qquad \PR-\text{almost surely}.
\]
\end{theorem}

\begin{remark}
 In the case considered in~\cite{BGP}, we may obtain the same results with $\alpha=\hatxi/\log \beta$ by the same methods. This exponent was conjectured in~\cite{BGP}. 
\end{remark}

\begin{remark} In dimension $d=2$, the  critical bias $\lambda_c:(p_c,1) \to (0,\infty)$ is monotone increasing in $p$.
\end{remark}

\subsubsection{Link to trap models} \label{weird_traps}
By this time, we feel that the reader should be familiar with the underlying mechanism that allows us to reduce a problem of directionally transient RWREs to the time spent in traps (we refer the reader to previous sections, see Section~\ref{trap_dep}, Section~\ref{sect_1dbtm} and Section~\ref{sect_arbre_trap} for more details). We recall that the key ingredients are that
\begin{enumerate}
\item the walk should not be able to go against its favored direction for long,
\item and that before advancing $n$ steps in the direction of the drift (time $\Delta_n$), the walk should visit a linear number of sites.
\end{enumerate}

Let us say a word about the actual techniques that had to be developed to formally link the biased random walk on a percolation cluster and sums of i.i.d.~random variables. The key is to justify the two previous properties and both can essentially be derived from proving that the biased random walk will exit large boxes in the direction of the drift (up to an error term exponentially small in the size of the box). This type of property has come to be known as condition $(T)$, $(T^{\alpha})$ or $(T')$ (see~\cite{SZ2} for a formal definition) and was introduced by Sznitman. He was motivated to introduce these conditions in light of a celebrated conjecture that any uniformly elliptic random walk in random environment that is directionally transient has positive speed, (whether or not the walk is reversible), see~\cite{Zeitouni} and~\cite{SZ2} for more on this conjecture.

This being said, we chose to discuss the geometry of traps. In the two-dimensional case, a trap is surrounded by a path of open dual edges, while, in higher dimensions, this surrounding is a surface comprised of open dual plaquettes. We refer to this as the trap surface. 

There is a significant difference in trap geometry between the cases $d =2$ and $d \geq 3$; see Figure~15. In any dimension, the typical trap is a long thin object oriented in some given direction. When $d\geq 3$, the trap surface is typically uniformly narrow, whereas in the two-dimensional case it is not.

\begin{figure}
\centering
\epsfig{file=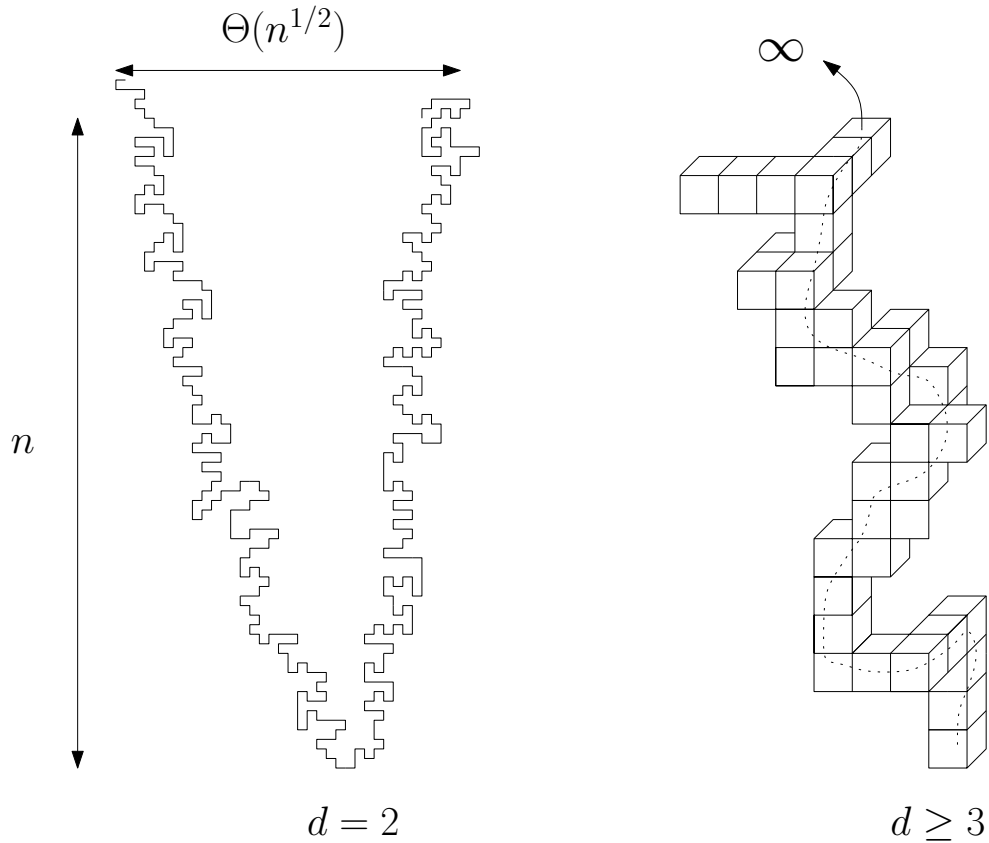, width=13cm}
\caption{The dual trap surfaces associated with typical traps.}
\end{figure}

This distinction means that, in the case that $d \geq 3$, there is a \lq\lq trap entrance\rq\rq, which is a single vertex, located very near the top of the trap, through which the walk must pass in order to fall into the trap. In contrast, in dimension $d=2$, the top of the trap is comprised of a line segment orthogonal to the trap direction and of length of the order of the square-root of the trap height. 

\begin{remark}
We wish to emphasize that for a generic non-axial $\vec{\ell}$, the trap direction does not coincide with $\vec{\ell}$. \end{remark}

\subsection{Open problems}

\subsubsection{The limiting speed}

Hardly anything is known about the limiting velocity, although qualitatively one would expect similarities with Figure 13. Here are some questions that one would like to see answered (in any dimension $\geq 2$).
\begin{question}
Is $v(\lambda,p)$ a unimodal function of $\lambda$ for every $p$? If so, for which $\lambda$ is the maximum achieved and what is the value of that maximum?
\end{question}

\begin{question}
Is $v(\lambda,p)$ an increasing (or non-decreasing) function of $p$ for every $\lambda$? For $p$ close to $1$, this has been known to be true, see~\cite{F} .
\end{question}

\begin{question}
What are the regularity properties of the function $v(\lambda,p)$?
\end{question}

Furthermore one would expect the Einstein relation to hold in the following sense: for every $p>p_c$
\[
\lim_{\lambda \to 0} \frac{v(\lambda,p)}{\lambda} =  C \sigma^2(p),
\]
where $C$ is a universal constant (independent of $p$) and $\sigma^2(p)$ is the variance appearing in the CLT  for the simple random walk on the supercritical percolation cluster of parameter $p$. The CLT for the simple random walk on the supercritical percolation cluster on $\Z^d$ is known to hold in the annealed sense, see~\cite{demasi}, and also in the quenched setting, see~\cite{sidosz} for $d\geq  4$ and the subsequent works of~\cite{BB} and~\cite{MP}.

\begin{question}
Does the Einstein relation hold for the biased random walk on a supercritical percolation cluster?
\end{question}

\subsubsection{Scaling limits}

We mention different conjectures.

\begin{conjecture} The CLT holds quenched for the biased random walk on the supercritical percolation cluster (see Theorem~\ref{the_theorem2}).
\end{conjecture}

The next conjecture is due to F.~and Alan Hammond.

\begin{conjecture} The right scaling for $X_n\cdot \vec{\ell}$ is $n^{\alpha}$ in the sense that $(X_n\cdot \vec{\ell}/n^{\alpha})_{n\geq 0}$ is tight. Moreover
\begin{itemize}
\item  If the bias direction is irrational (meaning that the line $\{\lambda \vec{\ell},\ \lambda \in\R\}$ in the unit torus is dense), the system should behave as a biased BTM, with stable scaling limits, aging ...
\item If the bias direction is rational (meaning that the line $\{\lambda \vec{\ell},\ \lambda \in\R\}$ in the unit torus is not dense), scaling limits (for the hitting time $\Delta_n$) in the zero-speed regime are expected to exist only on exponentially growing sub-sequences. 
\end{itemize}
A similar result should hold true in the regime of non-gaussian fluctuations.
\end{conjecture}

This conjecture can be understood in the light of Section~\ref{sect_arbre_lattice}. Indeed, the bias direction is irrational if, and only if, the height of a trap (i.e.~$\mathcal{B}\mathcal{K}$) has a distribution which is non-lattice. 

\begin{remark} Related models corresponding to `randomly biased random walks' in random environments on $\Z^d$ have been considered in the literature (see~\cite{Shen} and~\cite{fr} for a precise definition). Under certain conditions, a zero-speed regime related to trapping has been established (see~\cite{fr}), but the scaling limits of this model have not yet been studied.\end{remark}

\section{Random walks on critical trees.}\label{sec_crit}

In the previous sections, we have considered directionally transient (mainly biased) random walks on supercritical structures (Galton-Watson trees or percolation clusters). As we saw, dead-ends (or more generally poor connectivity) induced trapping, and the time spent by the walk in those traps had polynomial tails. By tuning the bias strength we could witness a rich range of behaviors starting from the traditional CLT to $\alpha$-stable scaling limits.

In critical structures, the environment is typically very  badly connected, in particular the dead-ends are typically orders of magnitude larger. This means that biased random walks on critical structures should have much more powerful slowdowns where actually the behavior of the walk is governed by the time spent in the largest trap. We will investigate this phenomenon in the next section.

As it turns out, when dealing with critical structures the natural setting for witnessing polynomial scaling (and a richer behavior) is to consider the simple random walk. Even though these models are not directionally transient we choose to discuss them, indeed they are very important models that are relevant to trapping but do not behave as classical (undirected) trap models. 

\subsection{Biased random walk on critical Galton-Watson trees}

We will present a result from~\cite{DFK} on the scaling limits for the biased random walk on critical Galton-Watson trees. 

\subsubsection{The model}\label{sect_crit_arbre}

We shall restrict ourselves to the case where $E[Z]=1$ and $E[Z^2]<\infty$ (excluding of course $Z= 1$ a.s.), whereas a slightly more general case is considered in~\cite{DFK}. 

In \cite{Kesten1}, Kesten showed that it is possible to make sense of conditioning a critical Galton-Watson tree to survive or `grow to infinity'. More specifically, there exists a unique (in law) random infinite rooted locally-finite graph tree $\mathcal{T}^*$, called the incipient infinite cluster (or I.I.C.) of the Galton-Watson tree, that satisfies, for any $n\in\mathbb{Z}_+$,
\[
E\left(\phi(\mathcal{T}^*|_n)\right)=\lim_{m\rightarrow\infty}E\left(\phi(\mathcal{T}|_n)|Z_{m+n}>0\right),
\]
where $\phi$ is a bounded function on finite rooted graph trees of $n$ generations, and $\mathcal{T}|_n$, $\mathcal{T}^*|_n$ are the first $n$ generations of $\mathcal{T}$, $\mathcal{T}^*$ respectively. We denote by ${\bf P}$ the law associated to  $\mathcal{T}^*$.

A key tool throughout this study is the spinal decomposition of $\mathcal{T}^*$ that appears as \cite[Lemma 2.2]{Kesten}, and which can be described as follows. First, $\mathbf{P}$-a.e. realization of $\mathcal{T}^*$ admits a unique non-intersecting infinite path starting at the root that we label $\rho_0,\rho_1,\rho_2,\ldots$. Conditionally on this `backbone', the number of children of vertices on the backbone are independent, each distributed as a size-biased random variable $\tilde{Z}$, which satisfies
\begin{equation}\label{sizebiasz}
\mathbf{P}\left(\tilde{Z}=k\right)=k\mathbf{P}(Z=k),\hspace{20pt} \text{ for all }k\geq 1.
\end{equation}
Moreover, conditional on the backbone and the number of children of each backbone element, the trees descending from the children of backbone vertices that are not on the backbone are independent copies of the original critical branching process $\mathcal{T}$. The critical trees are dead-ends and act as traps. See Figure~16 for a typical picture of the I.I.C.~of a Galton-Watson tree.

\begin{figure}
\centering 
\epsfig{file=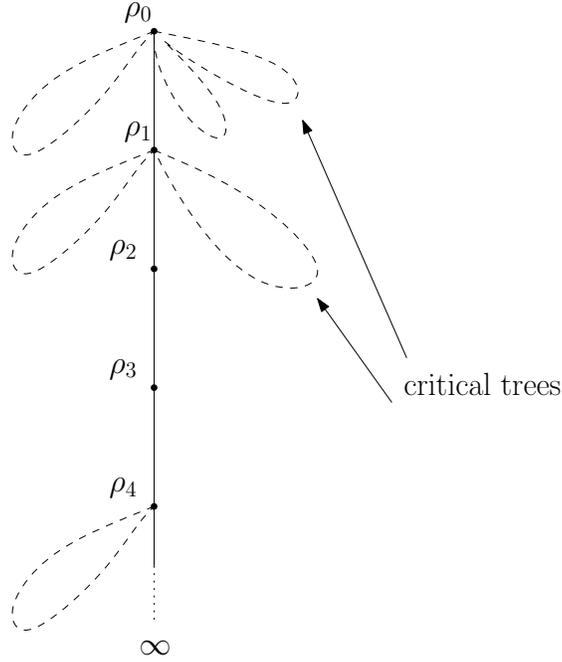, scale=0.7}
\caption{The I.I.C.~of a Galton-Watson tree. The dashed areas are critical Galton-Watson trees}
\end{figure}

On this random tree  $\mathcal{T}^*$, we can define $(X_n)_{n\geq 0}$, the $\beta$-biased random walk (for $\beta>1$) according to the rules stated in Section~\ref{sect_intro_biased}. We write $\PR$ for the associated annealed law.

\subsubsection{Results}

Let us define 
\[
\Delta_n:=\inf\left\{m\geq 0:\:X_m=\rho_n\right\}
\]
to be the first time the process $X$ reaches level $n$ on the backbone. The main result of~\cite{DFK} is the following functional limit theorem.

\begin{theorem}
Under the hypothesis mentioned above, we have that as $n\rightarrow\infty$, the laws of the processes
\[\left(\frac{\ln_+ \Delta_{nt}}{n\ln \beta}\right)_{t\geq 0}\]
under $\PR$ converge weakly with respect to the Skorohod $J_1$ topology on $D([0,\infty),\mathbb{R})$ to the law of an extremal process $(m(t))_{t\geq0}$.
\end{theorem}

We recall a way to define the law of an extremal process $m=(m(t))_{t\geq 0}$. Suppose that $(\xi(t))_{t\ge 0}$ is the symmetric Cauchy process, i.e., the L\'evy process with L\'evy measure given by $\mu((x,\infty))=x^{-1}$ for $x>0$,  we can then set
\[
m(t)=\max_{0<s\le t}\Delta\xi(s),
\]
where $\Delta\xi(s)=\xi(s)-\xi(s^-)$. (Observe that $(m(t))_{t\ge 0}$ is thus the maximum process of the Poisson point process with intensity measure $x^{-2}dxdt$.)

Another characteristic property that was shown in~\cite{DFK} is that the random walk also exhibits extremal aging (a concept that first appeared in~\cite{gun}).

\begin{theorem}
Under the hypothesis mentioned above, for any $0<a<b$, we have
\[
\lim_{n\to\infty}\PR\left(\pi(X_{e^{an}})=\pi(X_{e^{bn}})\right)=\frac{a}{b},
\]
where for any $t>0$, $\pi(X_t)$ denotes the projection of $X_t$ onto the backbone.
\end{theorem}

\begin{remark} It is interesting to note that in this context the lattice effect (see Section~\ref{sect_arbre_lattice}) cannot be felt. See Remark~\ref{rem_nolat} for more details.
 \end{remark}

\subsubsection{Links to trap models}

The decomposition of the critical tree described at the end of Section~\ref{sect_crit_arbre} shows that this model is essentially a biased random walk on $\Z$ with waiting times given by the excursion time $T_{\text{trap}}$ of a biased random walk in critical Galton-Watson trees.

Using the intuition from reversible Markov chain (more specifically the mean return formula), we will approximate $T_{\text{trap}} \approx \beta^H$, where $H$ is the height of the tree or, more formally, $H =\max \{n \geq 0, Z_n >0 \}$ where $Z_n$ is the size of the $n$-th generation. By a result from~\cite{KNS} we know that 
\[
{\bf P}[H\geq n]\sim \frac{ 2}{\text{Var}(Z)} n^{-1},
\]
which implies that ${\bf P}[T_{\text{trap}}\geq n]\sim C \ln^{-1} n$.

\begin{remark}\label{rem_nolat}  We emphasize that since $\ln(\cdot)$ is a slowly varying function it is true that $\ln n \sim \ln \lfloor n \rfloor$. This means that lattice effects (see Section~\ref{sect_arbre_lattice}) cannot be felt. Even more interestingly, we can see that the system becomes insensitive to a lot of information, indeed if $Y$ is a bounded random variable we can see that $Y\beta^H$ and $\beta^H$ have the same tail.
\end{remark}

This means that the $\beta$-biased random walk on a critical Galton-Watson tree behaves like a biased trap model with waiting times that have slowly-varying tails (of order $\ln^{-1}(n)$). Such models were studied in~\cite{DFK} and are shown to universally have extremal processes as scaling limits and exhibit extremal aging. In such types of models, at any large time the walk has an overwhelming probability to be in a large trap in which it has been for almost the entire time since it left the origin. This contrasts with the case of heavy-tails where the walk divides the majority of its time among a finite, but very large, number of traps.

\subsection{Simple random walk on critical Galton-Watson trees}

In this section we will consider a critical Galton-Watson tree conditioned to be infinite and perform a simple random walk on it. Actually, this is exactly the case $\beta=1$ described in Section~\ref{sect_crit_arbre}. To avoid repetition, we refer the reader to that section for the details about the model.

The behavior of the simple and the biased random walk on critical Galton-Watson trees are radically different. There are two main aspects to this difference
\begin{enumerate}
\item the time spent in the traps will have polynomially decaying tails instead of exponentially decaying ones,
\item the walk on the backbone is now recurrent which means that large traps will be revisited repeatedly.
\end{enumerate}

It was initially proved in~\cite{Kesten1} that a simple random walk $X_n$ on the I.I.C.~of a Galton-Watson tree (the infinite critical tree) is such that $\abs{X_n}/n^{1/3}$ converges in distribution to a non zero limit. However, the nature of the limiting distribution remained unknown for very long. The reason for this is that, although this model is a one dimensional model with traps, it cannot be approximated by a traditional undirected BTM. This can be seen since at a given site on the backbone, the random walker might, or might not, enter the trap deeply and hence the waiting time in a given trap may vary greatly  during different visits.

\begin{remark} Although it is also true, for biased random walks on supercritical (or critical) Galton-Watson trees, that the waiting time in a given trap may vary greatly  during different  visits, this does not influence the limit. The reason for this is the different nature of the embedded process, i.e.~the walk on the backbone. In the biased case, a lot of visits to the same trap may occur, but only in a short period of time after which the traps will never be visited again. In some sense, a big trap is only visited once. In the context of the simple random walk, the same trap may be revisited at much later times creating a more complicated behavior.
\end{remark}

The insufficiencies of the BTM to explain the limiting behavior of the simple random walk on a critical Galton-Watson tree prompted the introduction of a more general class of models called the Randomly Trapped Random Walks (RTRW). For these models the trapping landscape is not necessarily given by exponential random variables with random means (as in the BTM), but by randomly chosen random variables (i.e.~randomly chosen probability distributions).

 The possible scaling limits of such models are studied in~\cite{RTRW}. They obviously include the Fractional Kinetics (FK) models (see \cite{mula}) and the Fontes Isopi Newman (FIN) diffusion (see \cite{fin02}) that appear in the (non-directed) BTM (see~\cite{BenarousCerny}  for an overview of this topic). But those scaling limits also contain new processes called spatially subordinated brownian motions (SSBM) and, so-called, FK-SSBM mixtures. We refer to~\cite{ RTRW} for a precise definition of FK-SSBM.
 
 The limiting process of the simple random walk on the I.I.C.~of a critical Galton-Watson tree is (essentially) one of those SSBMs.

\subsubsection{Results}

We start by presenting the definitions necessary to describe the limiting process. Consider a $1/2$-stable subordinator $U_t$ and let $\rho$ be its random Lebesgue-Stieltjes measure on $\R$, i.e.~$\rho (a,b] = U_b-U_a$. As $U_t$ is pure jump process, $\rho$ is purely atomic and we can write
\[
\rho=\sum_{i\in \Z} \overline{y_i} \delta_{\overline{x_i}}.
\]

Then let $((S^i(t))_{t\geq 0})_{i\in \N}$ be an i.i.d.~ sequence of random processes having the (annealed) law of the inverse local time at the root of the Brownian motion on the Continuum random tree (see~\cite{croyd} and~\cite{aldous} for more precision on this object). We take this sequence to be independent of $\rho$.

Finally, let $(B_t^+)_{t\geq 0}$ be a one-dimensional standard Brownian motion reflected at the origin which is independent of the previously defined objects and let $l^+(x,t)$ be its local time. Then, define 
\[
\phi(t)=\sum_{i\in \N} \overline{y}_i^{3/2} S^i((1/2)\overline{y}_i^{-1/2}l^+(\overline{x}_i,t)).
\]

The main result of~\cite{bac} is that
\begin{theorem}
Let $X_n$ be a simple random walk  on the I.I.C.~of a Galton-Watson tree. For each $T>0$ we have that $(\epsilon X_{\lfloor \epsilon^{-3}t\rfloor})_{t\in [0,T]}$ converges in distribution as $\epsilon \to 0$ to $(B^+_{\phi^{-1}(t)})_{t\in [0,T]}$ on the space $D(0,T)$ equipped with the uniform topology. Here $\phi^{-1}$ denotes the generalized inverse of $\phi$, see below~(\ref{aboveq}).
\end{theorem}

\begin{remark} A similar result can be proved on the I.P.C.~also known as the invasion percolation cluster (see~\cite{chayes} for the first mathematical study of the I.P.C.) \end{remark}

\subsection{Simple random walk on the infinite incipient cluster of $\Z^d$}

For the reader who is unaware of the definition of percolation, we refer him to Section~\ref{perco}, where there is a quick presentation of the model.

\subsubsection{The model and the Alexander-Orbach conjecture}

In this section we will be discussing critical percolation, i.e.~the case $p=p_c(d)$. If $d = 2$ or $d \geq 19$, it is known (see for example~\cite{Grimmett}) that there is no infinite cluster at criticality and this is conjectured to hold for any $d \geq 2$. For large $n$, the local properties of these large finite clusters can, in certain
circumstances, be understood by looking  at them as subsets of an infinite cluster, called the incipient
infinite cluster (I.I.C.~for short). This is very much similar to the procedure described in Section~\ref{sect_crit_arbre}, where a similar construction is achieved on the tree. We do not wish to describe this procedure in further details and refer the reader to~\cite{Kesten} for the case $d=2$ and~\cite{Hoff} for $d\geq 19$. We denote the law of the I.I.C.~by ${\bf P}_{\text{IIC}}$.

We consider the simple random walk on the I.I.C.~of $\Z^d$ ($d\geq 2$).  Let us introduce the classical notion of spectral dimension of an infinite connected  graph $G$. It is defined to be
\[
d_s(G)=-2 \lim_{n\to \infty} \frac{\ln p_{2n}(x,x)}{\ln n}, \qquad \text{(if this limit exists)},
\]
where $x\in G$ and $p_n(x,x)=P^{G}[X_n=x\mid X_0=x]$ is the heat-kernel of the simple random walk on $G$. For a thorough overview of problems related to heat-kernel estimates we refer the reader to~\cite{Kumagai}.

One of the most famous physics  conjectures for the simple random walk on the I.I.C.~is known as the Alexander-Orbach conjecture~\cite{AO}. It states that for $d\geq 2$, the I.I.C.~of $\Z^d$ should have a spectral dimension of $4/3$. While
it is now thought that this is unlikely to be true for small $d$ (see~\cite{BJK} Section 7.4), it has been proved to hold for large dimensions in~\cite{KN}
\begin{theorem}\label{AO}
 For $d\geq 19$, the Alexander-Orbach conjecture holds true for the I.I.C.~of $\Z^d$, that is, ${\bf P}_{\text{IIC}}$-a.s.
\[
\lim_{n\to \infty} \frac{\ln p_{2n}(x,x)}{\ln n}=-\frac 23,
\]
furthermore denoting $\Delta_n=\inf\{i\geq 0,\ d_{\text{IIC}}(0,X_i)=n\}$, the first time that the walk hits a point at distance $n$ (in the I.I.C.~metric) from its starting point verifies
\[
\lim_{n\to \infty} \frac{\ln\Delta_n}{\ln n}=3.
\]
\end{theorem}

\begin{remark} The previous theorem indicates that the fluctuations of $X_n$ are of the order $n^{1/3}$. \end{remark}

\subsubsection{A sketch of proof}

Let us say a word about the proof. The following is highly non-rigorous but gives a good explanation for the exponent $2/3$. Let us assume that we are given a graph verifying
\begin{enumerate}
\item $R_{\text{eff}}(0,x_n)\approx n^{\rho}$ for any $x_n$ at distance $n$ of the origin. Here $R_{\text{eff}}(0,x)$ is the effective conductance between $0$ and $x$ (see~\cite{LP} or~\cite{DoyleSnell}). The I.I.C.~ in high dimensions is very thin and it is known that, in a certain sense, the previous relation holds for $\rho=1$.
\item $\text{Card} (B_{\text{IIC}}(0,n)) \approx n^{\nu}$, where $B_{\text{IIC}}(0,n)=\{x\in \text{I.I.C.},\ d_{\text{IIC}}(0,x)\leq n\}$. This is an information on the growth rate of large balls in the I.I.C.. In high dimensions, the I.I.C.~verifies a similar property for $\nu=2$.
\end{enumerate}

Let us now consider the commute time formula (see~\cite{commute}). On a finite graph $G$, for any two vertices $x$ and $y$, we have
\[
T(x\to y)+T(y \to x)=2R_{\text{eff}}(x,y) \text{Card} (E(G)),
\]
where $T(x\to y)$ denotes the expected hitting time of $y$ for the simple random walk starting from $x$ and $\text{Card}(E(G))$ the cardinal of the number of edges of $G$. Under reasonable assumptions, one would expect $T(x\to y)$ and $T(y\to x)$ to be of comparable orders. Leaving all rigor behind, the previous formula should  imply a formula similar to 
\[
T(0\to \partial B_{\text{IIC}}(0,n)) \approx R_{\text{eff}}(0, \partial B_{\text{IIC}}(0,n)) \text{Card} (B_{\text{IIC}}(0,n)) \approx n^{\rho+\nu},
\]
since $\text{Card} (E(B_{\text{IIC}}(0,n))) \approx \text{Card} (B_{\text{IIC}}(0,n))$. Note that this readily implies the second affirmation of Theorem~\ref{AO}. 

Let us explain how to obtain estimates on the heat-kernel from this. It takes time of the order of $n^{\rho+\nu}$ to reach $\partial B_{\text{IIC}}(0,n)$. Waiting twice the time it typically takes for the walk to reach $\partial B_{\text{IIC}}(0,n)$ (which is still of order $n^{\rho+\nu}$) we should leave enough time for the walk to mix inside $B_{\text{IIC}}(0,n)$ and be approximately uniformly distributed in that ball. Since $\text{Card} (B_{\text{IIC}}(0,n)) \approx n^{\nu}$, this means that $p_{n^{\rho+\nu}}(0,0) \approx n^{-\nu}$, or, stated differently
\[
p_{n}(0,0) \approx n^{-\nu/(\rho+\nu)}.
\]

Since for the I.I.C.~in high dimension we have $\nu=2$ and $\rho=1$, we have justified the result.

\subsection{Open problems}

The main open problems concerning random walks on critical structure are problems on the I.I.C.~of $\Z^d$. Let us mention the following problems.
\begin{question}
Prove that Alexander-Orbach does not hold for $\Z^2$ and find the correct scaling exponent.
\end{question}

\begin{question}
 For $\Z^d$ with $d\geq 19$, what are the scaling limits of the simple random walk on the I.I.C.?
 \end{question}
 
\begin{question}
What is the right scaling for the biased random walk on the I.I.C.~of $\Z^d$? 
\end{question}

Here, there are two possible definitions for a biased random walk. One could consider the cartesian bias, which has a definition similar to that of the biased random walk on the supercritical percolation cluster (see Section~\ref{perco}). Alternatively, one could consider a topological bias, where the walk experiences a bias along the chemical distance away from the origin. Both models should exhibit drastically different behaviors.

\section{Appendix on heavy-tailed random variables}

We give a quick overview of basic facts about sums of heavy-tailed i.i.d.~random variables. For more information on this subject we refer the reader to~\cite{Petrov} or~\cite{GK}.

When considering random walks in random environments and studying trapping we are naturally led to study sums of i.i.d.~random variables, i.e.~$S_n:=\sum_{i=1}^n X_i$ where $X_i$ is  a sequence of i.i.d.~random variables. In our context those random variables naturally happen to be non-negative so our discussion will be sometimes be restricted to that case.

The most basic properties for $S_n$ are the following
\begin{enumerate}
\item if $E[X]<\infty$, then $S_n/n \to E[X]$ almost surely by the law of large numbers,
\item if $E[X^2]<\infty$, then  $(S_n-nE[X])/n^{1/2}$ converges in distribution to a normal distribution $\mathcal{N}(0,\text{Var}{X})$ by the central limit theorem,
\item if $E[X^2]<\infty$, then the invariance principle implies that the process $(S_t^{(N)}, t\geq 0):=(N^{-1/2}(S_{\lfloor tN \rfloor}- tNE[X]), 0\leq t\leq T),$  converges as $N$ goes to infinity to  Brownian motions in the Skorohod topology, see Chapter 3 of~\cite{EK}.
\end{enumerate}

However, if $E[X^2]=\infty$, the previous results do not give any description for the fluctuations of $S_n$ and, if furthermore $E[X]=\infty$, even the law of large number fails. This case is commonly referred to as heavy-tailed.

 Our main goal here is to describe the limiting laws or processes that appear when trying to rescale $S_n$ when we do not have a finite first or second moment.

\subsection{Limiting laws}

Our main goal in the section is to understand conditions that can lead to the existence of constants $a_n$ and $b_n$ such that 
\begin{equation}\label{limstable}
\frac{S_n -b_n}{a_n}\xrightarrow{(d)} Y,
\end{equation}
where $Y$ is a non-degenerate random variable. Very much related to that question is a certain family of distributions called stable laws.

\subsubsection{Stable laws}\label{def_stable_law_section}

\begin{definition}
A non-degenerate distribution of a random variable $X$ is a stable distribution if, for any integer $k>0$ there exists constants $a_k$ and $b_k$ such that if $X_1,\ldots,X_k$ are i.i.d.~with the same distribution as $X$, then $(X_1+\ldots X_k-b_k)/a_k$ has the same distribution as $X$. \end{definition}

It turns out that stable laws is a four-parameter family. One of those parameters is usually very important, so that stable laws are usually referred to $\alpha$-stable laws, where $\alpha \in(0,2]$ is this central parameter.
\begin{definition} For $\alpha \in (0,2]$, an $\alpha$-stable random variable $\mathcal{S}_{\alpha}$ is defined by its characteristic function
\[
E[\exp({\bf i}t \mathcal{S}_{\alpha})]=\exp({\bf i}tc-b\abs{t}^{\alpha}(1+{\bf i} \kappa \text{sgn}(t) w_{\alpha}(t))),
\]
where $-1\leq \kappa \leq 1$, $b,c\in \R$ and $\text{sgn}(t)$  is the sign function and 
\[
w_{\alpha}(t)=\begin{cases} -\tan(\pi\alpha/2) & \text{ if } \alpha \neq 1 \\ (2/\pi) \log \abs{t} & \text{ if } \alpha=1 \end{cases}.
\]

When $\kappa=\pm 1$, we say that the stable law is completely asymmetric. 
\end{definition}

\begin{remark} The $2$-stable distribution is nothing more than the normal distribution with mean $c$ and variance $2b$.\end{remark}

\subsubsection{Limiting theorems}\label{sect_def_domain}

The importance of stable laws is related to the following (see~\cite{Durrett} Theorem 3.7.4).
\begin{theorem}
$Y$ is the limit in law of $(X_1+\cdots+X_k-b_k)/a_k$ for some i.i.d.~sequence $X_i$ if, and only if, $Y$ has a stable law.
\end{theorem}

In other words by re-centering and rescaling a sum of i.i.d.~random variables we can only obtain $\alpha$-stable limits, that is, in~(\ref{limstable}) the random variable $Y$ has to be stable.

\begin{definition} We say that a random variable $X$ is in the domain of attraction of an $\alpha$-stable law if there exist sequences $a_n$ and $b_n$ such that for some i.i.d.~sequence $X_i$ (with distribution $X$) we have that  $(X_1+\cdots+X_k-b_k)/a_k$  converges in distribution to an $\alpha$-stable law. Otherwise, we say the random variable $X$ does not belong to the domain of attraction of an $\alpha$-stable law.
\end{definition}

At this point, we know that only stable law can be obtained in the limit of~(\ref{limstable}), but, conversely, we would also like a sufficient criterion for the convergence of $S_n$.  For this let us make the following definition. 
\begin{definition}\label{def_sec_varlente}
 We say that a function $L(x)$ is slowly varying if $\lim_{x\to \infty} \frac{L(Kx)}{L(x)}= 1$ for any $K>0$.\end{definition}

\begin{remark} The function $L(x)=\ln x$ is slowly varying, whereas $x^{\epsilon}$ is not for any choice of $\epsilon>0$.\end{remark}

There exists an explicit criterion to verify if a random variable is in the domain of attraction of an $\alpha$-stable law (see~\cite{Durrett} Theorem 3.7.2)
\begin{theorem}\label{sumiidgen}
Suppose $X_1,X_2,\ldots$ are i.i.d.~with a distribution that satisfies 
\begin{enumerate}
\item $\lim_{x\to \infty} P[X_1>x]/P[\abs{X_1}>x]=\theta$,
\item $P[\abs{X_1}>x]=x^{-\alpha}L(x)$,
\end{enumerate}
where $\alpha<2$ and $L$ is slowly varying. Let $S_n=X_1+\cdots+X_n$,
\[
a_n=\inf \{x: P[\abs{X_1}>x]\leq n^{-1}\} \text{ and } b_n=nE[X_1\1{\abs{X_1}\leq a_n}].
\]

As $n\to \infty$, $(S_n-b_n)/a_n$ converges in distribution to an $\alpha$-stable law.
\end{theorem}

\begin{remark}\label{no_recenter} If $\alpha<1$, then we may chose $b_n=0$. \end{remark}

\begin{remark}\label{no_dom_attract} For $\alpha\in (0,2)$, if we cannot write $P[\abs{X}>x]=x^{-\alpha}L(x)$ (where $\alpha<2$ and $L$ is slowly varying), then $X$ does not belong to the domain of attraction of any stable law.
\end{remark}

We have now given a pretty detailed picture of what replaces the law of large numbers and the central limit theorem for sums of i.i.d.~random variables $X_i$ when the variance or the first moment of $X$ are infinite. We now turn to the description of the limiting processes that appear in the limit of sums of i.i.d.~random variables.

\subsection{Limiting processes}

Let us consider the two typical situations where $P[X\geq x]\sim x^{-\alpha}$ with $\alpha<1$ or $\alpha\in (1,2)$. Then, for any $t>0$
\[
\text{  if $\alpha<1$ then }\frac{S_{\lfloor tn\rfloor}}{n^{1/\alpha}} \xrightarrow{(d)} t^{1/\alpha}\mathcal{S}_{\alpha},
\]
 and
 \[
 \text{ if $\alpha \in (1,2)$, then }
\frac{S_{\lfloor t n\rfloor}-tnE[X]}{n^{1/\alpha}} \xrightarrow{(d)} t^{1/\alpha}\mathcal{S}_{\alpha} ,
\]
where $\mathcal{S}_{\alpha}$ is a certain $\alpha$-stable random variable.

It is a natural question to wonder how the whole process $(S_t^{(n)}, t\geq 0):=(n^{-1/\alpha}S_{\lfloor tn \rfloor}, 0\leq t\leq T)$, behaves as $n$ goes to infinity when $\alpha <1$. A similar question could be asked for $(S_t^{(n)}, t\geq 0):=(n^{-1/\alpha}(S_{\lfloor t n\rfloor}-tnE[X]), 0\leq t\leq T)$ when $\alpha \in (1,2)$.

In both cases, we have $S_t^{(n)}$ converges in the Skorohod topology (see Chapter 3 of~\cite{EK}) to an $\alpha$-stable L\'evy process.

\subsubsection{$\alpha$-stable L\'evy processes}\label{def_levyproc}

For a general treatment of L\'evy process, we refer the reader to~\cite{Bertoin}. The following definition will be sufficient for our purposes.
\begin{definition}
A c\`adl\`ag stochastic process $\{X_t,t\geq 0\}$ is called an $\alpha$-stable L\'evy process if 
\begin{enumerate}
\item $X_0=0$ a.s.,
\item $X$ has independent increments,
\item $X_t-X_s$ is distributed as $(t-s)^{1/\alpha} \mathcal{S}_{\alpha}$ where $\mathcal{S}_{\alpha}$ is an $\alpha$-stable distribution (see Section~\ref{def_stable_law_section} for the definition of stable laws).
\end{enumerate}
\end{definition}

\begin{remark} If the stable law appearing in the definition of an $\alpha$-stable L\'evy process is completely asymmetric we also say that the process is completely asymmetric. \end{remark}

In the case where we consider non-negative random variables that are in the domain of attraction of an $\alpha$-stable law (with $\alpha<1$), we know that the limiting distribution $S_n$ properly rescaled will be a positive $\alpha$-stable law (indeed for $\alpha<1$, we  do not recenter $S_n$, see Remark~\ref{no_recenter}). This means that the limiting $\alpha$-stable L\'evy process $S_t$ of $S_t^{(n)}$  is increasing. 

\begin{definition}
An increasing L\'evy process is called subordinator. Furthermore an increasing $\alpha$-stable L\'evy process is called $\alpha$-stable subordinator.
\end{definition}

These processes are of special interest to us since they occur naturally as limiting process when the trapping is strong enough to break the law of large numbers.

\subsubsection{Subordinators and aging}

 The class of subordinator can be described by a two-parameter family $(d,\mu)$ where $d\geq 0$ and $\mu$ is a measure on $(0,\infty)$ satisfying 
 \[
 \int_0^{\infty} (1\wedge x) \mu(dx)<\infty.
 \]
 
 The law of a subordinator  $S_t$ is then uniquely defined by its Laplace transform,
 \[
 E[\exp(-\lambda S_t)]=\exp(-t \Phi(\lambda)),
 \]
 where the Laplace exponent is 
 \[
 \Phi(\lambda)=d+\int_0^{\infty} (1-e^{-\lambda x})\mu(dx).
 \]
 
 All the subordinators appearing in these notes have no drift (i.e.~d=0) and are $\alpha$-stable (with $\alpha\in(0,1)$), which means that their Laplace exponent satisfies
 \[
 \Phi(\lambda)=c\lambda^{\alpha}=\frac{c\alpha}{\Gamma(1-\alpha)}\int_0^{\infty}(1-e^{-\lambda x})x^{-1-\alpha}dx,
 \]
 where $\Gamma$ is the usual Gamma-function.

An $\alpha$-stable subordinator without drift $S_t$ ($\alpha\in (0,1)$) has positive jumps. Let us assume that the Laplace exponent of $S_t$ is given by the previous equation. It is know that if one considers the random measure on $\R$ given by $\rho(a,b]=S_b-S_a$, then $\rho(dx)=\sum_i y_i\delta_{x_i}(dx)$ where $(x_i,y_i)$ is a Poisson point process  on $\R^+$ with intensity $\frac{c\alpha}{\Gamma(1-\alpha)} x^{-1-\alpha}dxdt$.

The jumps of an $\alpha$-stable subordinator appearing as a limit of a rescaled sum of i.i.d.~random variables represents the time spent on a certain trapping site. As such, assuming for simplicity that the scaling is $n^{1/\alpha}$, the probability that a walk is in the same trap at times $an^{1/\alpha}$ and $bn^{1/\alpha}$ ($a<b$) for large $n$ should be given by the probability that the limiting $\alpha$-stable subordinator jumps over the interval $[a,b]$. We are then able to link the aging property (see~(\ref{zd_5})) to the following key fact about stable subordinators.

\begin{theorem}
The probability that an $\alpha$-stable subordinator without drift $S_t$ jumps over the interval $[a,b]$ (i.e.~there is no $t\in \R$ such that $S_t \in [a,b]$) is equal to 
\[
P[S_{T(b)^-}<a]=P[\text{ASL}_{\alpha}\in [0,a/b]],
\]
where $T(b):=\inf\{t,\ S_t>b\}$ and $\text{ASL}_{\alpha}$ denotes the generalized arcsine distribution with parameter $\alpha$. The law $\text{ASL}_{\alpha}$  is supported on $[0,1]$ and its distribution function is given by
\[
P[\text{ASL}_{\alpha}\in [0,x]]= \frac{\sin{\alpha \pi}}{\pi} \int_0^{x} y^{\alpha-1}(1-y)^{-\alpha}dy.
\]
\end{theorem}

\section*{Acknowledgements}

The authors want to thank Benedikt Rehle for the simulations appearing in Figure 2-3-4. They are also indebted to Daniel Kious, Laurent Tournier and Olivier Zindy for many useful comments on an earlier version of these notes. Finally, the authors also wish to thank Elie A\"id\'ekon and Manuel Cabezas for interesting discussions.

\end{document}